\documentclass[reqno]{amsart}

\usepackage{hyperref}
\usepackage{amssymb}
\usepackage{dsfont}
\usepackage[all]{xy}

\usepackage{graphics}

\usepackage{epsfig}
\def\accentsfrancais{applemac}
   \RequirePackage[\accentsfrancais]{inputenc}

\newcommand\supp[1]{\ensuremath{{\rm supp}\left(#1\right)}}

\def\eps{\varepsilon}
\newcommand{\indiq}{{\bf 1}}

\def\R{\mathbb{R}}
\def\N{\mathbb{N}}
\def\Z{\mathbb{Z}}

\def\RR{\mathbb{R}}
\def\NN{\mathbb{N}}
\def\ZZ{\mathbb{Z}}

\newcommand{\cjk}{c^{i}_{j,k}}
\newcommand{\cla}{c_{\lambda}}
\newcommand{\pla}{\psi_{\lambda}}
\newcommand{\Laj}{\Lambda_j}
\newcommand{\dla}{d_{\lambda}}
\newcommand{\clap}{c_{\lambda '}}
\newcommand{\Hmin}{{h^{\Ome}_\EEE}}

\newcommand{\ep}{\epsilon}

\newcommand{\La}{\Lambda}

\newcommand{\EEE}{{\mathcal E}}
\newcommand{\ome}{\omega}

\newcommand{\Ome}{\Omega}
\newcommand{\bsp}{{\bf b}^{s, \infty}_p}
\newcommand{\btsp}{{\bf \tilde{b}}^{s, \infty}_p}
\newcommand{\Osc}{{\mathcal O}sc}
\newcommand{\al}{\alpha}

\newcommand{\la}{\lambda}
\newcommand{\lajx}{\la_j (x_0)}
\newcommand{\e}{\varepsilon}
\newcommand\sk{\smallskip}
\newcommand{\cJ}{{\mathcal J}}
\newcommand{\cL}{{\mathcal L}}
\newcommand{\cH}{{\mathcal H}}
\newcommand{\sm}{{s-}}
\newcommand{\as}{\mathrm{a.s.}}
\newcommand{\gen}[1]{\langle #1 \rangle}
\def\prob{\mathbb{P}}

\newtheorem{df}{Definition}
 \theoremstyle{plain}
\newtheorem{thm}{Theorem}
\newtheorem{prp}{Proposition}
\newtheorem{lem}{Lemma}
\newtheorem{cor}{Corollary}
\newtheorem{remark} {Remark}

\newsavebox{\fmbox}
\newenvironment{fmpage}[1]
 {\begin{lrbox}{\fmbox}\begin{minipage}{#1}}
 {\end{minipage}\end{lrbox}\fbox{\usebox{\fmbox}}}

\newcommand{\ho}{H\"older }
\newcommand{\BD}{\begin{df}}
\newcommand{\ED}{\end{df}}
\newcommand{\BT}{\begin{thm}}
\newcommand{\ET}{\end{thm}}
\newcommand{\BP}{\begin{prp}}
\newcommand{\EP}{\end{prp}}
\newcommand{\BL}{\begin{lem}}
\newcommand{\EL}{\end{lem}}
\newcommand{\BC}{\begin{cor}}
\newcommand{\EC}{\end{cor}}

\newcommand{\BE}{\begin{equation}} 
\newcommand{\EE}{\end{equation}}

\title{Local Multifractal Analysis}
\author{Julien Barral, Arnaud Durand,  St\'ephane Jaffard \and St\'ephane Seuret}
\date{July 30, 2012}
\address{LAGA,    Institut Galil\'ee, 
Universit\'e Paris 13 \\ 
99 avenue Jean-Baptiste Cl\'ement  \\ 
93430 - Villetaneuse \\  France} 
\email{barral@math.univ-paris13.fr}
\address{Laboratoire de Math\'ematiques d'Orsay, UMR 8628\\ Universit\'e Paris-Sud\\ 91405 Orsay Cedex\\ France} 
\email{arnaud.durand@math.u-psud.fr}
\address{Universit\'e Paris Est  \\ Laboratoire d'Analyse et de Math\'ematiques Appliqu\'ees, UMR 8050\\ 61 avenue du G\'en\'eral de Gaulle\\ 94010 Cr\'eteil Cedex\\ France}
\email{jaffard@u-pec.fr, seuret@u-pec.fr}

\begin{document}

\begin{abstract}
We introduce a local multifractal formalism adapted to functions, measures or distributions which display   multifractal characteristics that can change with time, or location. We  develop this formalism in a general framework   and  we work out several examples of measures and functions where this setting is relevant. 
\end{abstract}

\maketitle


\section{Introduction}

 Let $f$ denote a function, a positive Radon measure, or, more generally, a distribution defined on a nonempty open set $\Ome$.   One often  associates with $f$ a  pointwise exponent,  denoted by  $h_f  (x)$, which allows to quantify the  local smoothness of $f$ at $x$. On the mathematical side, the purpose of multifractal analysis is  to determine the fractal  dimension of the level sets of the function  $x\mapsto h_f  (x)$. 
  Let 
 \[ E_H = \{ x: \;\; h_f(x) =H\}. \]
  The {\em multifractal spectrum}  of $f$ (associated with  the regularity exponent $h_f$) is 
$$d_f (H) = \dim \, E_H  $$
(where $\dim$ denotes the Hausdorff dimension, see  Definition \ref{defmeshaus}). 
 Multifractal spectra yield a  description of the local singularities of the function, or measure, under consideration.

Regularity exponents (and therefore the multifractal spectrum) of many  functions, stochastic processes,  or measures used in modeling can be theoretically determined  directly from the definition.    However,  usually, one cannot recover these results numerically  on simulations,  because the  exponents  thus obtained  turn out to be    extremely erratic, everywhere discontinuous functions it is for instance the case of L\'evy processes \cite{Jaffard:1999fg}, or of multiplicative cascades (see the book \cite{BBBFHJMP}, and, in particular the review paper by J.~Barral,  A.~Fan,  and J. Peyri\`ere)  so that a direct determination of  $h_f (x)$ leads to totally instable computations. A fortiori, the estimation of the   multifractal spectrum from its definition  is unfeasable.  The multifractal formalism is a tentative way to bypass the intermediate step of the determination of the pointwise exponent, by relating the  multifractal spectrum directly with  {\em averaged quantities}  that are effectively computable on experimental data.   Such quantities can usually be interpreted as global regularity indices. For instance, the first one historically used in the function setting ($\zeta_f (p)$, referred to as {\em  Kolmogorov scaling function})  can be defined as follows; for the sake of simplicity, we only consider in this introduction  the function setting and  we assume here that the functions considered are defined on the whole $\RR^d$. 

 Recall that the Lipschitz  spaces  are defined, for 
  $s  \in (0,1)$, and  $p \in [1, \infty]$, by 
$ f \in \mbox{Lip}  ( s , L^p  (\RR^d)) $ if $f \in L^p$ and $\exists C>0$  such that $\forall h  >0$, 
\begin{equation} \label{nicol1}  \parallel f(x+h) -f(x)\parallel_{L^p} \leq C h^s .  \end{equation}
 (the definition for larger $s$ requires the use of higher order differences, and the extension to $p <1$ requires to replace  Lebesgue spaces by Hardy spaces, see \cite{Jaffard:2004fh}). 
Then 
\begin{equation} \label{nicol}   \zeta_f (p) =  p \cdot \sup \{ s: f \in \mbox{Lip}  ( s, L^p  (\RR^d)) \} .
\end{equation}

Initially introduced  by U. Frisch and G. Parisi in the mid 80s, the purpose of multifractal analysis is  to  investigate the relationships between  the pointwise regularity information  supplied by $d_f (H) $ and the global regularity information  supplied by $\zeta_f (p) $.  Note that these quantities can be computed on the whole domain of definition $ \Ome$  of $f$, or can be restricted to  an open subdomain $\ome \subset \Ome$. 
 A natural question is to understand how they depend on the  region $\ome$ where they are computed.  It is remarkable that, in many situations, there is no dependency at all on $\ome$; we will then say that the corresponding quantity is { \em homogeneous}.  It is the case for several classes of stochastic processes. For instance, sample paths of  L\'evy processes (and fields) \cite{Durand:2007fk,Durand:2010uq,Jaffard:1999fg}, L\'evy processes in multifractal time \cite{BSadv}, and  fractional Brownian motions (FBM) almost surely  have homogeneous H\"older spectra,    and, in the case of  FBM, the Legendre spectrum also is homogeneous, see \cite{JafJMP,BB2}.  In the random setting, it is also the case for many examples of multiplicative cascades, see \cite{BM04}. Many deterministic functions or measures  also are homogeneous (homogeneity is usually not explicitly stated as such in the corresponding papers, but is implicit in the determination of the spectra). This is for instance the case for self-similar or self-conformal measures when one assumes the so-called open set condition, or for Gibbs measures on conformal repellers (see for instance \cite{Olsen95,Pa,Pesin}). It is also the case for many applications, for instance  the Legendre spectra raising from natural experiments (such as turbulence, see \cite{AJW1,AADMV} and references therein)  are found to be homogeneous.
  
  On the opposite, many natural objects, either theoretical or coming from real data,  have been shown to be  non-homogeneous : Their multifractal characteristics  depends on the domain $\Ome$ over which they are observed:
\begin{itemize}

\item
It is the case of some classes of Markov processes, see \cite{Barral:2010vn} and Section \ref{sec:markov}, and also of some Markov cascades studied in \cite{Bae}. 

%
%
\item

Some self-similar measures when the open set condition is relaxed into the weak-separation condition may satisfy the multifractal formalism only when restricted to some intervals  (see \cite{Shm,FLW,Testud2,Fengmatrix2}).

\item
In applications, many types of signals, which have a human origin, can have multifractal characteristics that change with time: A typical example is supplied by finance data, see \cite{AJW1}, where changes can be attributed to  outside phenomena  such as political events,  but also to the  increasing sophistication of  financial tools, which may lead to instabilities (financial crises)  and implies that some characteristic features of the data, possibly  captured by multifractal analysis,  evolve with time.   This situation is also natural in image analysis because of the { \em occlusion phenomenon}; indeed, a natural  image is a patchwork of textures with different
characteristics, so that its  global spectrum of singularities  reflects the  multifractal nature of each component, and also
 of the boundaries (which may also   be fractal) where discontinuities appear.    Note that the notion of { \em local Hausdorff dimension }  which plays a central role in this section, has been introduced in \cite{JurSta} precisely with the motivation of image analysis.

\item Functions spaces with varying smoothness have been introduced motivated by the study of the relationship between general pseudo-differential operators and later by questions arising in  PDEs, see \cite{Schn} for a review on the subject; scaling functions  with characteristics depending on the location  are then the natural tool to measure  optimal regularity in this context. 
We will investigate this relationship in Section \ref{functional}. 

\end{itemize}

This paper will provide new examples   of  multifractal characteristics   which depend on the domain of observation. In such situations, the determination of a 
local spectrum  of singularities  for  each ``component''   $\ome \subset \Ome$ will carry more information than the knowledge of  the ``global'' one only.  A natural question is to understand how the different quantities  which we have introduced depend on the  region $\ome$ where they are computed.

 Some of the notions studied in this paper have been already introduced in \cite{Barral:2010vn}; 
let us  also mention that a  local $L^q$-spectrum was already introduced in \cite{KaRaSu}, where the authors studied this  notion for measures in doubling metric spaces (as well as the notion of local homogeneity) and obtained, for instance,   upper bounds for the dimensions of the sets of points with given lower and upper local dimensions using this local concepts. The goal of their approach was to investigate conical density and porosity questions. In our paper,  on top of measures, we also deal with functions,  get comparable upper bounds for the multifractal spectra, and the examples we develop are very different.

\medskip

 Let us now  make precise the notion we started with, namely pointwise regularity. 
 The two most widely used exponents are the pointwise H\"older exponent of functions and the local dimension of measures. In the following,  $B(x_0, r)$ denotes the open  ball of center $x_0$ and radius $r$.

\BD  Let $\mu$ be a  positive Radon  measure  defined on an open subset $\Ome \subset \R^d$.  
Let $x_0 \in \Ome$ and let ${\alpha} \geq 0$.   The  measure $\mu $   belongs to
$h^{\alpha}(x_0)$  if 
\BE  \label{HolMes}  
\exists C, R >0, \; \forall r \le R, \hspace{1cm}  \mu (B(x_0, r)) \leq C r^{\alpha}. \EE

Let $x_0$ belong to the support of $\mu$. The {lower local dimension } of $\mu$ at $x_0$ is 
 \BE \label{dmu1}
h_\mu(x_0)=\sup\{\alpha :\mu\in h^{\alpha}(x_0)\} = \liminf_{r\to 0^+} \frac{\log \mu(B(x_0,r))}{\log \, r}.
\EE
\ED

We now turn to the case  of  locally bounded functions.  In this setting, the notion corresponding to the lower local dimension is the pointwise  H\"older regularity. 

\BD
Let $x_0 \in \R^d$ and let ${\alpha} \geq 0$.  Let $f: \Ome \rightarrow \RR$ be a   locally bounded function; $f$  belongs to
$C^{\alpha}(x_0)$  if there exist
$C, R  >0$ and a polynomial
$P$  of degree at most $\alpha$  such that
\BE \mbox{ if }\;\;     | x-x_0| \leq R,  \;\; \mbox{ then }\;\; \hspace{1cm}  
\label{Hol} |f(x)-P(x-x_0)| \leq C |x-x_0|^{\alpha}.
\EE
The {H\"older exponent} of $f$ at $x_0$ is  
\BE
\label{defpointwise}
 h_f(x_0)=\sup\{\alpha:f\in C^{\alpha}(x_0)\}.
 \EE
\ED



\medskip


This paper is  organized as follows:

 In Section \ref{sec:locmulanal},  we recall the notions of dimensions that we will use (both in the global and local case),     we prove some basic results concerning the notion of local Hausdorff dimension,   and  we recall the  the wavelet characterization of pointwise H\"older regularity.

In Section \ref{sec:derMultForm}   we recall  the multifractal formalism on a domain  in a general abstract  form  which is adapted both to the function and the measure setting; then the  corresponding version of  {  local multifractal formalism} is obtained, and we draw its relationship with the notion of { \em germ  space}.

In Section \ref{sec:measures},   we investigate more precisely the local multifractal analysis of measures, providing natural and new examples where this notion indeed contains  more information than the single multifractal spectrum. In particular,  we introduce new cascade models  the local characteristics of which change  smoothly with the location; 
here again, we show that the local tools introduced in Section \ref{sec:locmulanal} yield the exact multifractal characteristics of these cascades.

 In Section \ref{locspecstoch}, we review the  results concerning some Markov processes which do not have stationary increments; then we show that the notion of local spectrum allows to recover the   exact pointwise behavior of the Multifractional Brownian Motion (in contradistinction with the usual ``global'' multifractal formalism).
 
In Section \ref{other}   we consider other regularity exponents characterized by dyadic families, and show  how they can be characterized in a similar way as the previous ones,  by log-log plot regressions of quantities  defined on the dyadic cubes.    


 Finally,  in Section \ref{functional}   the relationship between the local scaling function and   function spaces with varying smoothness is developed.


\section{Properties of the local Hausdorff dimension and the local multifractal spectrum}
  \label{sec:locmulanal}

  \subsection{Some notations and recalls}

In order to make precise  the different notions of multifractal spectra, 
we need to   recall the notion of dimension which will be used.

\BD  \label{defmeshaus} Let   $A\subset \RR^d$. 
If $\ep>0$  and  $\delta \in [0,d]$,  we denote 
\[ M^\delta_{\ep} =\inf_R \;  \left( \sum_{  i} | A_i |^\delta \right) ,\]
where $R$ is an  $\ep$-covering of   $A$, i.e. a covering of  $A$  by bounded sets  $\{ A_i\}_{i \in
\NN}$
 of diameters  $| A_i | \leq \ep$.
The infimum is therefore taken on all  $\ep$-coverings $R$.

For any $\delta \in [0,d]$, the 
$\delta$-dimensional  Hausdorff measure of 
$A$ is 
\[  mes_\delta (A) = \displaystyle\lim_{\ep\rightarrow 0}  M^\delta_{\ep}. \]
There exists
 $\delta_0 \in [0,d]$ such that 
\[ \begin{array}{l} \forall \ 0< \delta < \delta_0, \;\;
 mes_\delta (A) = + \infty  \;\; \;\;  \mbox{ and} \;\; \;\;\;\;   \forall \delta > \delta_0, \;\;
 mes_\delta (A) = 0
; \end{array}\]
this critical  value $\delta_0$ is called the  Hausdorff dimension of  $A$, and is  denoted by $\dim (A)$. 
By convention, we set  $\dim (\emptyset) = -\infty$.

\ED  

In practice, obtaining lower bounds for the Hausdorff dimension directly from the definition involve considering all possible coverings of the set, and is therefore not practical.  One rather  uses the  {\em mass distribution principle} which 
  involves instead the construction of a well-adapted  measure.
  
\begin{prp} Let $A \subset \RR^d$ and let 
 $\mu$ be a Radon measure such that $\mu (A) >0$; 
  \[  \mbox{ if}\;\; \forall x \in A, \;\;\; \limsup_{ r \rightarrow0} \frac{\mu (B(x,r))}{r^s} \leq C \;\;  \;\mbox{ then }\;\; \;\; { \mathcal H}^s (A) \geq \frac{\mu (A)}{c}.\] 
\end{prp}
  We will see in Section \ref{sec:locmulspec} a local version of this result.

  \medskip

 Apart from the Hausdorff dimension, we will also need another notion of dimension: The   {\em  packing dimension} which was introduced by  C. Tricot, see \cite{Tric1}:

\BD \label{defpack}
 Let   $A$ be a bounded subset of $ \RR^d$; if 
 $\ep>0$,  we denote by  $N_\ep (A)$   the smallest number of sets of radius $\ep$ required to cover $A$. 
The {  lower box dimension}  of $A$ is 
\[ \underline{\dim}_B (A)  =\liminf_{\ep \rightarrow 0} \frac{\log N_\ep (A)}{-\log \ep}  .\]
 The {packing dimension} of a set $A\subset \RR^d $  is 
\BE \dim_p (A) = \inf \left\{ \sup_{i\in \NN}  \left( \underline{\dim}_B A_i: A \subset \bigcup_{i=1}^\infty A_i \right)\right\} \EE
(the infimum is taken over all possible partitions of $A$ into a countable collection $A_i$). 
\ED

  \subsection{Local Hausdorff dimension  }

\


 In situations where the spectra are not homogeneous, the purpose of multifractal analysis is to understand how they  change with the location where they are considered. In the case of the multiractal spectrum, this amounts to  determine how the Hausdorff dimension of the set $E_f(H)$  changes locally. This can be performed using the   notion of { \em local Hausdorff dimension}, which can be traced back  to \cite{JurSta} (see also \cite{Bae} where this notion is shown to be fitted to the study of  { \em deranged Cantor sets}). 
 
\BD Let $A \subset \RR^d$, and $x \in \RR^d$. The local Hausdorff dimension of $A$ at $x$ is 
  the function defined by 
\BE \label{dimloc}  
  \forall x \in  \overline{A}, \hspace{1cm} 
\dim (A, x) = \lim_{r\rightarrow 0}\;  \dim (A \cap {B (x,r)}). 
\EE
\ED 

{ \bf Remarks:} 
\begin{itemize}
\item   The limit exists because, if  $ \Ome_1 \subset \Ome_2$, then   $ \dim ({\Ome_1}) \leq  \dim ({\Ome_2})$;  therefore the right-hand side of (\ref{dimloc}), being  a non-negative increasing function of $r$,   has a limit when $r \rightarrow 0$. 
\item  We can also conisder this quantity as defined on the whole $\RR^d$, in which case, it takes the value $-\infty$ outside of $\overline{A}$.
\item The same  definition allows to define a  local dimension, associated with any other definition of fractional dimension; one gets for instance a notion of { \em local  packing dimension}.  
\end{itemize}
   
 The following result shows that the local Hausdorff dimension  encapsulates all the information concerning the Hausdorff dimensions of the sets of the form $A \cap \omega$, for any open set $\ome$. 
 
 \BP \label{dimloc1}
Let $A \subset \RR^d$;
then for any open set $\ome$ which intersects $A$,
\BE \label{dimloc2}Ê  
  \dim (A \cap \omega )  = \sup_{ x\in \ome} \; \dim (A, x) .\EE
\EP

\begin{proof}
For $r$ small enough, $B_r \subset \omega$; it follows that   
\[ \forall x \in \omega, \hspace{1cm}  \dim (A, x) \leq  \dim (A \cap \omega ), \]
and therefore 
 $\displaystyle\sup_{ x\in \ome} \dim (A, x) \leq  \dim (A \cap \omega )$. 

Let us now prove the converse inequality.  Let $(K_n)_{n \in \NN}$ be an increasing sequence of compact sets such that $\cup K_n = \ome$; then 
\[  \dim (A \cap \omega ) = \lim_{n \rightarrow \infty} \dim (A \cap K_n ). \] 
Let $\delta >0$ be given;  then
\[  \forall x \in K_n, \; \;   \exists r (x) >0, \hspace{1cm}   \dim (A \cap {B (x,r)}) - \dim (A, x) \leq \delta . \]
We extract a finite covering of $K_n$ from the collection $ \left\{ B(x, r(x) ) \right\}_{ x \in K_n}$ which yields a finite number of points
$x_1, \cdots x_N \in \omega$ such that  $ K_n \subset \bigcup  B(x_i, r(x_i) ) $; thus 
\begin{eqnarray*}
 \dim (A \cap K_n ) &   \leq  & \sup_{i=1, \cdots , N}  \dim (A \cap  B(x_i, r(x_i) )  )   \\
& \leq &\sup_{i=1, \cdots , N} \dim (A, x_i) + \delta  \leq \sup_{ x\in \ome} \dim (A, x) + \delta .
 \end{eqnarray*}
Taking $\delta \rightarrow 0$ and $N \rightarrow \infty$ yields the required estimate.  
\end{proof}

Proposition \ref{dimloc1} implies  the following  regularity for the local Hausdorff dimension.

\BC
Let $A$ be a given subset of $\RR^d$; then the function $x \rightarrow \dim (A,x) $ is upper semi-continuous.
\EC

\begin{proof}We have 
\[   \dim (A,x) = \lim_{r\rightarrow 0}\;  \dim (A \cap {B (x,r)}) = \displaystyle\lim_{r\rightarrow 0}\;  \sup_{ y\in B(x, r ) } \dim (A, y) = \limsup_{y \rightarrow x} \; \dim (A, y).  \]
\end{proof}

\subsection{Wavelets and wavelet  leaders}

\label{obowb}


In Section \ref{sec:derMultForm}  we will describe a general framework for deriving a multifractal formalism adapted to pointwise regularity exponents. The key property of these exponents  that we will need is that  they are derived from log-log plot regressions of quantities defined on the dyadic cubes. Let us first check  that it is the case for the pointwise exponent of measures.

  Recall that a  dyadic cube of scale  $j \in \Z $ is   of the form 
\BE \label{defla} \la   = \left[ \frac{k_1}{2^j}, \frac{k_1+1}{2^j}\right) \times \dots \times \left[ \frac{k_d}{2^j}, \frac{k_d+1}{2^j}\right),
\EE  where $k=(k_1, \dots k_d )\in \Z^d$. 
Each point $x_0 \in \R^d $ is contained in  a unique dyadic cube of scale $j$, denoted by $\lajx$. 

 Let $3\la_j(x_0) $ denote   the  cube with the same center as   $\la_j(x_0)$ and three times wider;  it is easy to check that  \eqref{HolMes}  and \eqref{dmu1} can be rewritten as
 $$  h_\mu(x_0) = \liminf_{j \rightarrow + \infty}  \frac{ \log 
  \mu (3\la_{j}(x_0) )  }{\log  \, 2^{-j}  }  .$$

We now show that  the H\"older exponent of a function can be recovered in a similar way, from quantities derived from wavelet coefficients. 
  Recall that orthonormal wavelet bases on $\RR^d$ are of the following form: There exist a function  $\varphi $ and  $2^d-1$ functions $\psi^{(i)}$ with the following properties: The  $\varphi (x-k)$ ($k\in \ZZ^d$)  and  the $2^{dj/2}\psi^{(i)}(2^jx-k)$
($k\in
\ZZ^d,$  $j\in \ZZ$) form an orthonormal basis of $L^2 (\RR^d)$. 
This basis  is $r$-smooth if $\varphi$ and the $\psi^{(i)}$ are  $C^{r}$ and if the  $\partial
^{\al}\varphi$, and the  $\partial
^{\al}\varphi  \psi^{(i)} $,  for  $| \al | \leq r$, have fast decay.
 Therefore, $\forall f \in L^2$, 
\begin{equation} \label{ecrit2} f(x) =  \sum_{k\in \ZZ^d}  c_k \varphi (x-k) +  \sum_{j=0}^{\infty} \sum_{k\in \ZZ^d}
\sum_{i} \cjk \psi^{(i)} (2^jx-k); \end{equation} the $c_k$ and $\cjk$ are the { \em wavelet coefficients } of 
$f$:  
\begin{equation} \label{cjk}  \cjk = 2^{dj} \int_{\RR^d} f(x)  \psi^{(i)} (2^j x -k) dx , \hspace{6mm} \mbox{ and} \hspace{6mm}
 c_k = \int_{\RR^d} f(x)  \varphi ( x -k) dx .\end{equation}

Note that  (\ref{cjk})  makes sense even if $f$ does not belong to $L^2$; indeed,  when using  smooth enough wavelets, these formulas 
can be interpreted as a duality product between  smooth functions (the wavelets) and distributions.

   Instead of   the three indices $(i,j,k)$,   wavelets will be indexed by  { dyadic cubes}  as follows:    Since the wavelet index 
$i$ takes 
$2^d-1$ values, we can assume that it takes values in  $\{ 0, 1 \}^d -(0, \dots , 0)$; we will use the notations
\[ \la\; (= \la (i,j,k)) \; = \displaystyle\frac{k}{2^j} + \displaystyle\frac{i}{2^{j+1}} 
 + \left[ 0, \displaystyle\frac{1}{2^{j+1}}\right)^d, \;\; 
\cla = \cjk, \;\; \pla (x)= \psi^{(i)}(2^jx-k) .
\]  Note that the cube $\la$ which indexes the wavelet   gives information about its location and scale; if one uses compactly supported wavelets, then $\exists C >0$ such that $\supp\pla \subset C \cdot \la$.

Finally,  $\Laj$ will denote the set of dyadic cubes  $\la$  which index a wavelet of scale  $j$, i.e. wavelets 
 of the form  $\pla (x)= \psi^{(i)}(2^jx-k) $  (note that  $\Laj$ is a subset of the dyadic cubes of side 
 $2^{j+1}$). 
We take for norm on $\RR^d$
\[ \mbox{if}\;\;  x = (x_1, \dots , x_d), \;\;\; \; |x| = \sup_{ i=1, \dots , d} | x_i |; \]
so that the diameter of a dyadic cube of side  $2^{-j}$ is  exactly $2^{-j}$. 

 In the following, when dealing with H\"older regularity of functions, we will always assume that, if a function $f$ is  defined on an unbounded set $\Ome$, then it has slow increase, i.e. it satisfies
 \[ \exists C, N >0  \hspace{1cm} |  f(x) | \leq C (1+ | x|)^N ; \]
   and, if $\Ome \neq \RR^d$, then  the wavelet basis used is compactly supported, so that, if $x_0 \in \Ome$, then the wavelet coefficients ``close''  to $x_0$ are well defined for $j$ large enough.


 Let $f$ be a locally bounded function, with slow increase.  The pointwise H\"older regularity of $f$ is characterized in terms of the
 { \em  wavelet leaders } of $f$:  
\begin{equation}\label{eq:defWL}
\dla = \sup_{ \la '  \subset  3 \la}  |\clap | .
\end{equation}

The assumptions we made on $f$ imply  that 
wavelet leaders are  well defined and  finite.  

 We note  $ d_j (x_0) =  d_{ \la_j (x_0)} $. The following result  allows to characterize the H\"older exponent  by the decay rate of the  $ d_{ \la_j (x_0)}$  when $j \rightarrow +
\infty$, see  \cite{Jaffard:2004fh}.

\BP \label{regponc} Let  $\al >0$ and let $\psi_\la$ be an orthonormal basis  with regularity $r > \al$. 
If    there exists $ \ep >0$ such that
$f\in C^\ep (\Ome)$,  then  
 \BE \label{leaderregpon} \forall x_0  , \hspace{1cm}  h_f (x_0) = \liminf_{j \rightarrow + \infty}  \frac{ \log 
  d_{ \la_j (x_0)}   }{\log  \, 2^{-j}  }  .
\end{equation}
 \EP
 
Hence, the pointwise \ho  exponent can be computed from a dyadic family. This is also the case for the lower dimension of a measure $\mu$. Indeed, if $3\la_j(x_0) $ stands for the dyadic cube $\la$ of generation $j$ neighboring $\la_j(x_0)$,  it is easy to check that  \eqref{HolMes}  and \eqref{dmu1} can be rewritten as
 $$  h_\mu(x_0) = \liminf_{j \rightarrow + \infty}  \frac{ \log 
  \mu (3\la_{j}(x_0) )  }{\log  \, 2^{-j}  }  .$$

\section{ A local multifractal formalism for a dyadic family}\label{sec:derMultForm}

\subsection{Multifractal analysis on a domain $\Omega$}

%

\BD  Let $\Ome$ be a non-empty open  subset of $\RR^d$.  A  collection  of nonnegative quantities $\EEE = (e_\la)$    indexed  by  the set of  dyadic cubes  $\la \subset \Ome$ is called a dyadic function on $\Ome$. 
\ED
 
 The choice of the dyadic setting may seem arbitrary; however, it is justified by two reasons:
 \begin{itemize}
\item  It is the natural choice when dealing with orthonormal   wavelet bases (though wavelets could be defined using other division rules,  in practice the dyadic one is the   standard   choice), and also the measure setting.
\item  When analyzing experimental data  through regressions on log-log plots, for a given resolution, the dyadic splitting yields the largest number of scales available in order to perform the regression. 
\end{itemize}

%
%
%

  \BD  \label{defexpabstra} The  pointwise exponents associated with a dyadic function $\EEE$  on $\Ome$ are the function $ {h} (x)$ and   $ \tilde{h} (x)$ : $ \Ome \rightarrow \RR$ defined  for  $x \in \Ome$  as follows: 
\begin{itemize}
\item The   { lower  exponent  }  of  $\EEE$ is
\BE \label{exploclimsix} h_\EEE (x) =  \liminf_{j \rightarrow + \infty}  \frac{ \log 
   e_{\la_j (x)}   }{\log  \, 2^{-j}  }  \end{equation}
\item  The    { upper   exponent }   of  $\EEE$   is 
\BE \label{exploclimsixbis}
   \tilde{h}_\EEE  (x) =      \limsup_{j \rightarrow + \infty}  \frac{ \log 
  e_{\la_j (x)}  }{\log  \, 2^{-j}  }   .
\end{equation}
\end{itemize}
By convention one sets $h_\EEE(x) = \tilde{h}_\EEE  (x) =+\infty $ if $x\notin$Supp($\EEE) $.

\ED

%
%
%
%
%
%

We saw in the introduction the first example of scaling function which has been used. We now define them in the abstract setting supplied by dyadic functions. 
 We denote by $\La_j^\Ome$ the subset of $\La_j$ composed of the  dyadic cubes contained in $\Ome$.

 \BD \label{defscalfuncloc}  Let $\Ome$ be a nonempty bounded open subset of $\RR^d$.
 The   structure function of a dyadic function $\EEE$  on $\Ome$  is defined  by 
\BE \label{defstrucfunc}  \forall p \in \RR ,  \hspace{1cm} S_j (\Ome, p) =   {\sum_{ \la \in \La_j^\Ome}} (e_\la)^p.  \EE

The  scaling function  of $\EEE$ on $\Ome$ is defined by 
\BE \label{scalfunc}   \forall p \in \RR , \hspace{1cm}  \tau_{\EEE}^\Ome (p)  = \liminf_{j \rightarrow + \infty}  \frac{ \log  S_j (\Ome, p)  }{\log\, 2^{-j} }      .
\end{equation}
 \ED
 
 
 
  If $\Ome$ is not bounded, one  defines the scaling function as follows: 
 \BE \label{omenotboun} \mbox{ if} \;\;  \Ome_n = \Ome \cap B(0, n), \;\;    \forall p \in \RR , \hspace{1cm}  \tau_{\EEE}^\Ome (p)  = \lim_{n \rightarrow  \infty}   \tau_{\EEE}^{\Ome_n} (p). \EE
 Note that the limit exists because the sequence is decreasing. From now on, we will   assume that  {\bf the set $\Ome$ is bounded}, so that, at each scale $j$, a finite number only of  dyadic cubes $\la$ satisfy $\la \subset \Ome$.  The corresponding results when  $\Ome$ is unbounded   follow easily from (\ref{omenotboun}).  \\Ê

Apart from the scaling function, an  additional ``global'' parameter  plays an important role for  classification in many applications; and, for multifractal analysis, checking  its positivity is a prerequisite  in the wavelet setting   (see \cite{AJW1} and references therein): The {\em uniform regularity exponent} of $\EEE$ is defined by 
\BE 
\label{caracbeswav3hol}  
\Hmin =\liminf_{j \rightarrow + \infty} \;\;  
 \;\;  \frac{ \log (
  \displaystyle  \sup_{\la \in \La_j }   e_\la) }{\log \, 2^{-j} }. 
\EE

The scaling  function $\tau_\EEE^\Omega$    is concave  (as a liminf of concave functions) taking values in $\overline{\RR}  $.   
  The following regularity assumption is often met in practice, and implies that $\eta^\Ome_{\EEE}$ is finite for any value of $p$. 
 
 \BD A dyadic function $\EEE$ is regular in $\Ome$  if 
  \BE \label{defab}  \exists C_1, C_2 >0, \; \exists A, B \in \RR\:  \hspace{3mm} \forall \la \subset \Ome : \;\; e_\la \neq 0 \Longrightarrow  \; \; \; 
C_1 2^{-Bj} \leq e_\la \leq C_2 2^{-Aj}. \EE 
 \ED 
 
 The existence of $A$ is equivalent to the condition $\Hmin  > -\infty$. More precisely, 
 \[ \Hmin = \sup\{ A: \; \mbox{ the right hand side of }  \;  (\ref{defab}) \; \mbox{ holds} \}  . \]
 
%
%
%
\medskip

In the measure case and in the H\"older exponent case, one can pick $A =0$. In the H\"older case, the uniform regularity assumption means that $A >0$. 
When the $e_\la$ are wavelet leaders, the assumption on the lower bound implies that the function $f$ considered has no $C^\infty$ components.  

\medskip

  Since the scaling function is concave, there is no loss of information in rather considering its  Legendre transform, defined by 
  \begin{equation} \label{formult1.1}
{ \mathcal L}^\Ome_\EEE (H) := \inf_{p \in \RR} ( Hp - \tau_{\EEE}^\Ome (p)  ).
\end{equation} 
The function ${\mathcal L}_\EEE^\Ome (H) $ is called the {\em Legendre spectrum} of $\EEE$.

Though it is mathematically equivalent to consider ${\mathcal L}_\EEE^\Ome (H) $ or $\tau_{\EEE}^\Ome (p)$, one often prefers to  work with the Legendre spectrum, because of its interpretation in terms of regularity exponents supplied by the {\em multifractal formalism}.

\BD \label{specsing}   Let $\EEE $ be a  dyadic function on $\Ome$, and define, for $H\in [-\infty, +\infty]$,  the level set associated with $\EEE$
\[ E^\Omega_\EEE(H) = \{ x \in \Omega: \hspace{3mm}  h_\EEE (x) =H\} . \]
$$ d_\EEE^\Ome: H \in \RR \mapsto  \dim  \ E^\Omega_\EEE(H). $$
 \ED
%

Let us now show how a heuristic  relationship  can be  drawn   between   the multifractal   and the Legendre spectra.
The definition of the scaling function (\ref{scalfunc})
 roughly means that, for  $j$ large,
 $  S_j (\Ome, p) \sim 2^{- \tau_{\EEE}^\Ome (p) j}.  $
Let us  estimate the contribution to $ S_j (\Ome, p) $ of the dyadic cubes  $  \la$ that cover the points of $E_\EEE(H)$.
By definition of $E_\EEE(H)$, they satisfy
$ e_\la  \sim 2^{-Hj}; $
 by definition of $d_\EEE^\Ome (H)$,   since we use cubes of the same width $2^{-j}$  to cover $E$, we need about $2^{ d_\EEE^\Ome (H) j}$ such cubes; therefore the corresponding
contribution   is  
$ \sim  2^{d_\EEE^\Ome(H) j}2^{-Hpj}$  $ = 2^{-j(Hp- d_\EEE^\Ome (H) )}.$ 
 When $j \rightarrow +
\infty$,  the smallest exponent brings an exponentially  dominant contribution, so that 
\begin{equation} \label{etamuleg} \tau_{\EEE}^\Ome (p)= 
\inf_H ( Hp - d_\EEE^\Ome (H)).\end{equation}
This formula can be interpreted as stating that the scaling function is the Legendre transform of the spectrum. 
Assuming that  $d_\EEE^\Ome (H)$ is concave, it can be recovered by an inverse Legendre transform, leading to   
\begin{equation} \label{formult1}
d_\EEE^\Ome (H) = \inf_{p \in \RR} ( Hp  - \tau_{\EEE}^\Ome (p)  ).
\end{equation} 
When this equality holds,  the dyadic function $\EEE $ { \em satisfies the multifractal formalism}, which therefore  
 amounts to state  that the Legendre spectrum coincides with the multifractal spectrum.

Note that the  derivation  we sketched  is not a mathematical proof, and the determination of the  range of validity  of
(\ref{formult1}) (and of its variants) is one of  the main mathematical problems concerning  multifractal analysis.  The only results which hold in all generality are upper bounds of dimensions of singularities. 
%
%
%
%
%
%
%
%
%
%

\BP \label{majspecprop}\cite{BMP,Jaffard:2004fh, JLVV} Let $\EEE$ be a dyadic function on $\Ome$.  Then
\begin{equation} \label{theo112} d^{\Ome}_\EEE ( H) \leq { \mathcal L}^\Ome_\EEE (H).
\end{equation}
 \EP

An important consequence of this corollary is supplied by the only case  where the knowledge of the scaling function  is sufficient to deduce the multifractal spectrum,  and even the pointwise exponent $h_\EEE$ everywhere.

\BC \label{contrex}   Let  $\EEE$ be a dyadic function. If its   scaling function    $\tau^\Ome_\EEE$  satisfies
\BE \label{valid} \exists \al >0\hspace{3mm} \mbox{ such that} \hspace{3mm}  \forall p \in \RR,  \hspace{3mm} \tau^\Ome_\EEE (p) = \tau_\EEE(0)  + \al p,\EE
then  the multifractal formalism is satisfied on $\Ome$, and the  lower exponent of $\EEE$  
 is 
\[   \forall x \in \mbox{Supp} \ \EEE ,  \hspace{1cm} h_\EEE (x)  = \al .  \]
\EC

\begin{proof} (of Corollary \ref{contrex})
 Assume that \eqref{valid} is true.
Then  $L_\EEE (H) = -\infty$ except for $H=\al$;  Corollary \ref{majspecprop} implies in this case that $d_\EEE(H) \leq -\infty$ for $H \neq \al$. Therefore only one H\"older exponent is present, so that $\forall x, \; h(x) = \al$;   it follows that $ d^{\Ome}_\EEE (\al) =1$, and the multifractal formalism  therefore holds. 
\end{proof}

This corollary has  direct implications in modeling: Indeed, several experimental signals  have   a linear scaling function.  In such situations,  multifractal analysis yields that the data have  a constant  pointwise exponent; therefore it   supplies  a non-parametric method  which allows to conclude  that modeling by, say, a fractional Brownian motion,  is appropriate (and the slope of the scaling function supplies the index of the FBM), see e.g.  \cite{AJW1} where one example of internet traffic data is shown.  We will also see a local version of Corollary \ref{contrex} which has implications in modeling: Corollary \ref{contrex2}.

 \subsection{Local multifractal formalism}  \label{sec:locmulspec}

%

\BD  \label{defpoint2}  Let $\EEE$ be a dyadic function on $\Ome$. The local multifractal spectrum of $\EEE$  is  the function defined by 
\BE \label{equivlocspec}  
\forall H,\;  \forall x \in \Ome,  \hspace{1cm} 
d_\EEE (x, H) = \dim (E_\EEE(H), x)  \left(  = \lim_{r\rightarrow 0} d^{B (x,r)}_\EEE ( H) \right) . 
\EE
\ED

The following result, which is a direct consequence of Proposition \ref{dimloc1}, shows that the local spectrum   
allows  to recover  the  spectrum of all possible restrictions of $\EEE$
on a subset $\omega  \in \Ome$.

\BC \label{spectresup}
Let $\EEE$ be a dyadic function on $\Ome$.
Then for any open set $\ome \subset \Ome$,
\BE \label{spectresup1}Ê \forall H\in \RR,  
\hspace{1cm}   d^{\omega}_\EEE ( H) = \sup_{ x\in \ome} d_\EEE (x, H) .\EE
\EC



 \BD
 A dyadic family $\EEE$  is said to be homogenously multifractal when   the local multifractal spectrum $d_\EEE(x,\cdot)$ does not depend on $x$, i.e.
 \[  \forall x \in \Ome , \; \forall H \in \RR, \qquad  d_\EEE (x, H) = d^{\Omega}_\EEE ( H) . \]
 \ED
 
%
  
   A local scaling function can also be defined  by making the set $\Ome$ shrink down to  the point $x_0$. 
   
   \BD
    \label{def:locscalfunc}  Let $\EEE$ be a dyadic function on $\Ome$. The local scaling  of $\EEE$  is  the function defined by 
\BE \label{equ:locscalfunc}  
\forall H,\;  \forall x \in \Ome,  \hspace{1cm} 
\tau_\EEE (x, p) = \lim_{r\rightarrow 0}   \tau^{B (x,r)}_\EEE (p) . 
\EE
   \ED 

Note that the right-hand side of (\ref{equ:locscalfunc}) is a decreasing function of $r$, and therefore it  has a limit when $r \rightarrow 0$.  Similarly as in the multifractal spectrum case, a  straightforward compacity argument yields that the  scaling function on any domain $\omega$ can be recovered from the local scaling function.

\BC  
Let $\EEE$ be a dyadic function on $\Ome$.
Then for any open set $\ome \subset \Ome$,
\BE \label{spectresup2}Ê \forall H\in \RR,  
\hspace{1cm}   \tau^{\omega}_\EEE (p) = \inf_{ x\in \ome} \tau_\EEE (x, p) .\EE
\EC

 \BD
 The  scaling function of a dyadic family $\EEE$  is said to be homogenous when   the local scaling function $\tau_\EEE(x,\cdot)$ does not depend on $x$. 
 \ED

 The upper bound supplied by Corollary  \ref{majspecprop}  holds for any  given ball $B(x, r)$. Fixing $x\in \Omega $ and making $r \rightarrow 0$, we obtain a following  local version of this result: 
\begin{equation} \label{locmulfor} \forall x \in \Ome, \; \forall H, \hspace{1cm} d_\EEE (x, H) \leq \inf_{ p\in \RR} \left( Hp -\tau_\EEE (x, p)    \right) . \end{equation} 
We will say that the { \em multifractal formalism holds locally} at $x$ whenever (\ref{locmulfor}) is an equality.

As above, this result has an important consequence: In some cases, it allows to determine the regularity exponent at every point, even in situations where this exponent is not constant.

\BC \label{contrex2}   Let  $\EEE$ be a dyadic function. If there exists a 
function $  \al :\RR \mapsto \RR$ such that  the   local scaling function    $\tau_\EEE$   satisfies  
\BE \label{valid2}     \hspace{3mm}  \forall x\in \Ome , \;  \forall p \in \RR,  \hspace{3mm} \tau_\EEE (x, p) =  \tau_\EEE(x,0)+ \al (x) p,\EE
then  the multifractal formalism is locally satisfied on $\Ome$, and the  lower exponent of $\EEE$ is 
\BE \label{muform2}  \forall x \in \Ome ,  \hspace{1cm} h_\EEE (x)  = \al (x) .  \EE
\EC   

This result is  a direct  consequence of (\ref{locmulfor}) and Corollary \ref{contrex}: Indeed, if (\ref{valid2}) holds, then   (\ref{locmulfor}) implies that  $d_\EEE (x, H)= -\infty $  if $H \neq \al (x)$.  We pick now an $H \neq \al (x)$;  recall that 
$ d_\EEE (x, H) = \lim_{r\rightarrow 0} d^{B (x,r)}_\EEE ( H) $;   
therefore $\exists R >0$ such that $\forall  r  \leq r$, $d^{B (x,r)}_\EEE ( H) = -\infty$. In particular, $H$ is not the pointwise  exponent at $x$. Since this argument holds for any $H \neq \al (x)$, (\ref{muform2}) holds, and Corollary \ref{contrex2}  follows.   \\ 

 We will see an application of Corollary \ref{contrex2}  concerning the multifractional Brownian Motion in Section \ref{sec:mbm}. 
Combining (\ref{locmulfor})   with Proposition \ref{spectresup}, yields the following upper bound.  

\BC \label{majspecprop2} Let $\EEE$ be a dyadic function on $\Ome$; for any open set $\ome \subset \Ome$,
\begin{equation} \label{theo113} \forall H,  \hspace{1cm} d^{\omega}_\EEE ( H) \leq \sup_{ x\in \ome}\;  \inf_{ p\in \RR} \left( Hp -\tau_\EEE (x, p) \right) . \end{equation} \EC  

It is remarkable that, though this result is a consequence  of Corollary \ref{majspecprop}, it usually yields a sharper bound.  Indeed, assume  for example  that the multifractal formalism  holds for two separated regions $\omega_1$ and $\omega_2$ yielding two different spectra $d_1 (H)$ and $d_2 (H)$; then 
(\ref{theo113}) yields $\max ( d_1 (H), d_2 (H))$ whereas the global multifractal formalism applied to $\Ome =  \omega_1 \cup \omega_2$ only yields  the concave hull of $\max ( d_1 (H), d_2 (H))$.  More generally, each time (\ref{theo113}) yields a non-concave upper bound, it will be strictly sharper than 
the result supplied by Corollary \ref{majspecprop}. 

\sk\sk

The uniform regularity  exponent also has a local form: 

\BD
The local exponent associated with $\EEE$ is the function 
\[ h_\EEE (x) = \lim_{r \rightarrow 0}  h_\EEE^{B(x,r)}. \]
\ED 

Note that the most general possible local exponents are lower semi-continuous functions, see \cite{JLVS}. 
 
It would be interesting to obtain a similar characterization for the functions $(x, H) \rightarrow d_\EEE (x, H)$ and $(x, p) \rightarrow  \tau_\EEE (x, p) $ (considered as as functions of two variables)
and determine their most general form. 

\subsection{An example from ergodic theory}\label{exampleergodic}
Let $\Omega=(0,1)$. Consider a $1$-periodic functions $\phi:\R\to \R$, as well as two  continuous functions $\gamma: [0,1] \to (0,\infty)$ and $\theta:[0,1]\to\R$. Let $T: x\in\R\mapsto 2x$. For $x\in\R$ and $j\in\N$ denote by $S_j\phi (x)$ the $j^{\text{th}}$ Birkhoff sum of $\phi$ at $x$, i.e., 
$$
S_j\phi(x)=\sum_{k=0}^{j-1}\phi(T^k x).
$$
Then, for any dyadic subinterval $\lambda$ of $\Omega$ of generation $j$, let
$$
e_\lambda= \sup_{x\in \lambda}e^{-\gamma (x)S_j\phi(x)-j\theta (x)}.
$$
When the functions $\gamma$ and $\theta$ are constant, the multifractal analysis of the dyadic family $\mathcal E=(e_\lambda)_{\lambda\subset \Omega}$ reduces to that  of the Birkhoff averages of $\gamma\phi+\theta$, since $\liminf_{j\to\infty} \frac{\log e_{\lambda_j(x)}}{\log 2^{-j}}=H$ if and only if $\liminf_{j\to\infty}S_j(x)/j=  \frac{H\log (2)-\theta}{\gamma}$. This is a now classical problem in ergodic theory of hyperbolic dynamical systems, which is well expressed through the thermodynamic formalism. The function $\log (2) \tau_{\mathcal E}^\Omega$ is the opposite of the pressure function of $-(\gamma\phi+\theta)$, that we denote by $P_{\gamma,\theta}(q)$, i.e.
\begin{eqnarray*}
-\log (2) \tau_{\mathcal E}^\Omega(p)=P_{\gamma,\theta}(p)&=&\lim_{j\to\infty}\frac{1}{j}\log\sum_{\lambda\in \Lambda^\Omega_j} \big (\sup_{x\in \lambda}e^{-\gamma S_j\phi(x)-j\theta}\big )^p\quad (p\in\R),\\
&=& P(-\gamma p) -\theta p,
\end{eqnarray*}
where $P=P_{-1,0}$ is the pressure function of $\phi$; and the following result follows for instance from \cite{FF}.

\begin{thm}
Let $H\in\R$; then $E^{\Omega}_{\mathcal E}(H)\neq\emptyset$ if and only if $H$ belongs to the interval $[{(\tau_{\mathcal E}^\Omega)}'(\infty),{(\tau_{\mathcal E}^\Omega)}'(-\infty)]$ and in this case $\tau_{\mathcal E}^\Omega(H)=\inf\{Hp-\tau_{\mathcal E}^\Omega(p):p\in \R\}$. 
\end{thm}
Continuing to assume that $\gamma$ and $\theta$ are constant, and using the fact that $\mathcal{E}$ possesses the same almost multiplicative properties as weak Gibbs measures (see \cite{K,FO} for the multifractal analysis of these objects), i.e. some self-similarity property, it is easily seen that we also have $\tau_{\mathcal E}^\omega=\tau_{\mathcal E}^\Omega$ and $d_{\mathcal E}^\omega=d_{\mathcal E}^\Omega$ for all open subsets of $\Omega$. 

Now suppose that $\gamma$ or $\theta$ is not constant. Such a situation should be seen locally as a small perturbation of the case where these functions are constant, and it is indeed rather easy using the continuity of $\gamma$ and $\theta$ to get the following fact.
\begin{prp} $\displaystyle \forall\ x\in\Omega,\ \forall q\in\R$,
\begin{equation}\label{taux}
 \tau_{\mathcal E}(x,p)=-\frac{P_{\gamma(x),\theta(x)}(p)}{\log (2)}= \frac{-P(-\gamma(x)p )+\theta(x) p}{\log (2)} .
 \end{equation}
\end{prp}

Suppose also that $\phi$ is not cohomologous to a constant,  i.e. the pressure function $P$ of $\phi$ is not affine, which is also equivalent to saying that the interval $I=[P'(-\infty),P'(\infty)]$ of possible values for $\liminf_{j\to\infty}S_j(y)/j$,  is non trivial. 

For all $H\in \R$, define 
$$
\xi_H:y\in (0,1)\mapsto \frac{H\log (2)-\theta(y)}{\gamma(y)}.
$$
Notice that $\liminf_{j\to\infty} \frac{\log e_{\lambda_j(y)}}{\log 2^{-j}}=H$ if and only if $\liminf_{j\to\infty}S_j(y)/j=h$ and $H=(\gamma(y)h +\theta (y))/\log(2)$, i.e. $h=\xi_H(y)$. 

Now fix $x\in (0,1)$. For $r>0$  we thus have  
\begin{equation}\label{low}
E^{B(x,r)}_{\mathcal E}(H)=\{y\in B(x,r): \liminf_{j\to\infty}S_j(y)/j=\xi_H(y)\}, 
\end{equation}
and due to Theorem 2.3 in \cite{BQ}, for all $H>0$, 
$$
\dim E^{B(x,r)}_{\mathcal E}(H)\ge \sup\{\inf\{P(p)-p\alpha:p\in\R\}: \alpha\in \mathrm{rg}({\xi_H}_{|B(x,r)})\cap \mathrm{int}(I)\}.
$$ 
Fix $H\in (\tau_{\mathcal E}'(x,\infty), \tau_{\mathcal E}'(x,-\infty))=(\gamma (x)P'(-\infty)+\theta(x),\gamma (x)P'(\infty)+\theta(x))$. By construction,  
\[{\xi_H}_{|B(x,r)}(x)=(H\log(2)-\theta(x))/\gamma(x)\in \mathrm{rg}({\xi_H}_{|B(x,r)})\cap \mathrm{int}(I). \]
 Thus, due to \eqref{low}, 
 \[ \dim E^{B(x,r)}_{\mathcal E}(H)\ge \inf\big \{P(p)-p\big (H\log(2)-\theta(x))/\gamma(x)\big ):p\in\R\big \}, \]
  which, due to \eqref{taux}, is exactly $\inf\{Hp- \tau_{\mathcal E}(x,p):p\in\R\}$. Since this estimate holds for all $r>0$, 
  \[
d_{\mathcal E}(x,H)\ge \inf\{Hp- \tau_{\mathcal E}(x,p):p\in\R\}, \] 
hence, by \eqref{locmulfor}, it follows that 
\[  
d_{\mathcal E}(x,H)=\inf\{Hp- \tau_{\mathcal E}(x,p):p\in\R\}.\]

For the case where  $H\in \{\tau_{\mathcal E}'(x,\infty), \tau_{\mathcal E}'(x,-\infty)\}$, it is difficult to conclude in full generality. We thus have proved the following result.
\BT
Suppose $\phi$ is not cohomologous to a constant. Fix $x\in \Omega$ and  $H\in\R$. If $H\not\in[\tau_{\mathcal E}'(x,\infty), \tau_{\mathcal E}'(x,-\infty)]=[\gamma (x)P'(-\infty)+\theta(x),\gamma (x)P'(\infty)+\theta(x)]$  then $E^{B(x,r)}_{\mathcal E}(H)=\emptyset$ for $r$ small enough, and if $H\in (\tau_{\mathcal E}'(x,\infty), \tau_{\mathcal E}'(x,-\infty))$ then $d_{\mathcal E}(x,H)=\inf\{Hp- \tau_{\mathcal E}(x,p):p\in\R\}$. 
\ET
Let us mention that if the union of the sets of  discontinuity points of $\gamma$ and $\theta$ has Hausdorff dimension 0, then the study achieved in \cite{BQ} shows that the previous result holds at any point $x$ which is a point of continuity of both $\gamma$ and $\xi$. Also, when $\phi$ and $\theta$ are positive, the family $\mathcal E$ can be used to build wavelet series whose local multifractal structure is the same as that of $\mathcal E$.

%
%
%
%
%
%
%
%
%

\section{Measures with varying local spectrum}
\label{sec:measures}

\subsection{General considerations}

Let $\mu$ be a positive Borel measure supported by $[0,1]^d$. recall that one derives from $\mu$  the dyadic family $\EEE_\mu= \{ e_\lambda:=\mu( 3\la)\}_{\la\in \Lambda}$.

It is obvious that  the definition \eqref{dmu1} of the local dimension $h_\mu(x_0) $ is equivalent to   \eqref {exploclimsix} with the dyadic family $\EEE_\mu$.  Similarly, the classical formalism for measures on $[0,1]^d$ is the same as the one described in the previous section for the family $\EEE_\mu$ on $\Omega=[0,1]^d$. Hence one can define a local multifractal spectrum for measures by Definition \ref{equivlocspec}.

In the measure setting, the following  result  shows that the mass distribution principle has a local version.

\BP
Let $\mu$ be a Radon measure, $A\subset \RR^d$ and  $x \in \overline{A} \cap \supp\mu$. Then 
\[ \dim (x, A) \geq h_\mu (x) . \]
\EP

\begin{proof}
It follows from (\ref{regunifmes}) applied on $A \cup B(x, r)$,  remarking that the hypothesis $x \in \supp\mu$ implies that $\mu (A \cup B(x, r)) >0$ and then letting $r \rightarrow 0$. 
\end{proof}

\medskip



We introduced the local multifractal spectrum to study non-homogeneous multifractal measures. It is interesting to recall the result of \cite{BucSeu}, where it is proved that homogeneous multifractal measures and non-homogeneous multifractal measures do not exhibit the same multifractal properties.
 
\begin{thm}
\label{thdarboux}
Consider a non-atomic homogeneous multifractal  measure supported on $[0,1]$. Then the intersection of the support of the (homogeneous) multifractal spectrum of $d_\mu$ with the interval $[0,1]$ is necessarily an interval of  the form $(\alpha,1] $ or $[\alpha,1] $, where $0\leq \alpha \leq 1$.
 \end{thm}

This is absolutely not the case for non-homogenouely multifractal measures: consider for instance two uniform Cantor sets $C_0$  and $C_1 $ of dimension 1/2 and 1/4 on the intervals $[0,1/2)$ and $[1/2,1]$. Then the barycenter of the two uniform measures naturally associated with  $C_0$ and $C$ satisfies 
\[ d_\mu(h) = \left\{ \begin{array}{rl}  
1/4 & \mbox{ if } \; h = 1/4, \\    1/2 &  \mbox{ if } \;  h = 1/2, \\   -\infty  &  \mbox{else.} \end{array}\right. \]
Hence the local spectrum is {\bf the} natural tool to study non-homogeneous multifractal measures. 

\subsection{A natural example where the notion of  local spectrum is relevant}

The Bernoulli (binomial) measure is perhaps the most natural and simple multifractal object, and it is now folklore that is is homogeneously multifractal. We make a very natural modification in its construction, which will break homogeneity by making the Bernoulli parameter $p$ depend on the interval which is split  in the construction. Doing this, we obtain a "localized" Bernoulli measure whose local spectrum  depends on $x$. This example is closely related with the example developed in Section~\ref{exampleergodic}.

\sk\sk

Let $p=[0,1] \mapsto (0,1/2)$ be a continuous mapping. For  $n\geq 1$, $({\ep_1,\ep_2, ...,\ep_n}) \in \{0,1\}$, we denote the dyadic number $k_{\ep_1\ep_2...\ep_n} =   \sum_{i=1}^n \ep_i 2^{-i}$ and  the dyadic interval $I_{\ep_1\ep_2 ...\ep_n} =\left[  k_{\ep_1\ep_2...\ep_n} ,k_{\ep_1\ep_2...\ep_n} +2^{-n}\right)$, where $n\geq 1$, $({\ep_1,\ep_2, ...,\ep_n}) \in \{0,1\}$, and we will use the natural tree structure of these intervals using the words $({\ep_1 \ep_2  ... \ep_n})$.  

 Consider the sequence of measures $(\mu_n)_{n\geq 1}$ built as follows:

\begin{itemize}
\item
$\mu_1$ is uniformly distributed on $I_0$ and $I_1$, and $\mu_1(I_0) = p(2^{-1})$ and $\mu_1(I_1) = 1-p(2^{-1})$.

\item
$\mu_2$ is uniformly distributed on the dyadic intervals $I_{\ep_1\ep_2}$ of second generation, and 
$$\mu_2(I_{\ep_10 }) =\mu_1(I_{\ep_1}) \cdot  p(k_{\ep_11}) \ \mbox{ and } \  \mu_2(I_{\ep_1 1 }) =\mu_1(I_{\ep_1}) \cdot   (1-p(k_{\ep_11})) .$$

\item
...

\item
$\mu_n$ is uniformly distributed on the dyadic intervals $I_{\ep_1\ep_2...\ep_n}$ of   generation $n$, and 
\begin{eqnarray*}
  \mu_n(I_{\ep_1\ep_2   ... \ep_{n-1} 0 })  & = &  \mu_{n-1}(I_{I \ep_1\ep_2   ... \ep_{n-1}}) \cdot p(k_{I \ep_1\ep_2   ... \ep_{n-1} 1})  \\
 \mbox{ and } \  \mu_n(I_{  \ep_1\ep_2   ... \ep_{n-1} 1 })  & =   & \mu_{n-1}(I_{\ep_1})  \cdot  (1-p(k_{ \ep_1\ep_2   ... \ep_{n-1}1} )). 
\end{eqnarray*}  
\end{itemize}

Observe that by construction, for every $n$, for every $p\geq n$ and every dyadic interval $I$ of generation $n$, one has $\mu_p(I)=\mu_n(I)$. 

\BD
The sequence of measures $(\mu_n)_{n\geq 1}$ converges weakly to  a measure $\mu$ that we call the "localized" Bernoulli measure associated with the map $p$.
\ED

Obviously, if $p$ is constant, one recovers the usual Bernoulli measure with parameter $p$.

We indicate the sketch of the proof to obtain the local multifractal properties of $\mu$. We do not use exactly the exponent $h_\mu$ defined by \eqref{dmu1}, for simplicity we work with the {\em dyadic} local exponent defined by 
$$h^d_\mu(x) =\liminf _{j\to +\infty} \frac { \log \mu(I_j(x)}{   \log \, 2^{-j}},$$
where (as usual) $ I_j(x)$ stands for the unique dyadic interval of generation $j$ containing $x$. What we are going to prove  also holds for the exponent $h_\mu$, but would require long technical developments. In particular, we would need an extension of  Corollary 2 of \cite{BQ} on localized multifractal analysis of Gibbs measures. This exponent $h^d_\mu$ can also be encompassed in the frame of Section \ref{sec:derMultForm} by using the dyadic family $\EEE=\{\mu(\la)\}_{\la\in \La}$, thus all the "local" notions we introduced hold for this exponent.

\BT
\label{theolocalmeas}
For every $x\in [0,1]$, the local spectrum associated with the exponent $h^d_\mu$ of $\mu$ at $x $ is that of a Bernoulli measure of a parameter $p(x)$, i.e.
$$\forall \ H\geq 0, \ \ \ \ d_\mu(x,H) = d_{\mu_{p(x)}}(H).$$
\ET

\sk\sk

For every $x \in [0,1]$, we consider its dyadic decomposition $x=\ep_1\ep_2 .... \ep_n ...$,  $\ep_i\in\{0,1\}$. Let $N_{0,n}(x) = \#\{1\leq k \leq n:  \ep_k=0\}$ and $N_{1,n}(x) = \#\{1\leq k \leq n:  \ep_k=1\} \ (= n-N_{0,n}(x) )$. We consider the asymptotic frequencies of 0's and 1's in the dyadic decomposition of $x$ defined as
$$ N_0(x) = \limsup_{n\to +\infty} \frac1 n N_{0,n}(x).$$

\BP
\label{prop:locmeasure}
For every $x \in [0,1]$, we have $$h^d_\mu(x) = -  N_0(x) \log_2 p(x)  -  (1-N_0(x)) \log_2 (1-p(x)).$$
\EP
Essentially, the localized binomial measure looks locally around $x$ like the binomial measure of parameter $p(x)$. 
\begin{proof}
Let us fix $q\in (0,1/2)$, and consider the classical Bernoulli measure $\mu_q$ of parameter $q$ on the whole interval $[0,1]$.
It is classical that the \ho exponent of $\mu_q$ at every point $x$ is 
\BE
\label{exphmuq}
h^d_{\mu_q}(x) = -  N_0(x) \log_2 q  -  (1-N_0(x)) \log_2 (1-q).
\EE
Inspired by this formula, a Caesaro argument gives the proposition. Indeed, by construction, the value of the $\mu$-mass of the interval $I_n(x)$ is given by
\begin{eqnarray*}
\mu(I_n(x)) & = &    \prod_{i=1}^{n}  p(k_{\ep_1\ep_2...\ep_{i-1}1})^*,
\end{eqnarray*}
where  
$$p(k_{\ep_1\ep_2...\ep_{i-1} 1})^* = \begin{cases} p(k_{\ep_1\ep_2...\ep_{i-1}1}) & \mbox { if } \ep_i=0\\
 1-p(k_{\ep_1\ep_2...\ep_{i-1}1}) & \mbox { if } \ep_i=1
\end{cases}.$$
Hence, 
\begin{eqnarray*}
\mu(I_n(x)) & = &    2^{   \sum_{i=1: \ep_i= 0}^n \log_2 p(k_{\ep_1\ep_2...\ep_{i-1}1}) + \sum_{i=1: \ep_i=0}^n  \log_2(1- p(k_{\ep_1\ep_2...\ep_{i-1}1}))},
\end{eqnarray*}

Since $p(k_{\ep_1\ep_2...\ep_{i-1}1}) $ tends to $p(x) $ when $i$ tends to infinity, and since $N_{0}(x)$ is the asymptotic frequency of zeros in the dyadic expansion of $x$, one sees that 
$$ \frac1n \sum_{i=1: \ep_i= 0}^n \log_2 p(k_{\ep_1\ep_2...\ep_{i-1}1}) \longrightarrow_{n\to +\infty}   N_0(x)p(x).$$
Similarly, since $p(x)<1/2$, 
$$ \frac1n \sum_{i=1: \ep_i= 1}^n \log_2 (1-p(k_{\ep_1\ep_2...\ep_{i-1}1})) \longrightarrow_{n\to +\infty}    (1- N_0(x))(1-p(x)).$$

Let $\alpha=-  N_0(x) \log_2 p(x)  -  (1-N_0(x)) \log_2 (1-p(x))$. The latter proves that, given $\ep>0$, there exists an integer $N$ such that $n\geq N$ implies that
$$2^{-n(\alpha +\ep)} \leq \mu(I_n(x)) \leq 2^{-n(\alpha-\ep)}.$$
This yields the result.
\end{proof}

Consider an interval $J\subset [0,1]$, and the multifractal spectrum $d_\mu(H,J) = \dim \, \{x\in J: h^d_\mu(x)=H\}$. The value of this spectrum is a consequence of  the following theorem of Barral and Qu  in \cite {BQ} (who proved this result for any Gibbs measure $\mu$).

\BT
\label{th:BQ} 
Fix $q\in (0,1/2)$, and consider the Bernoulli measure with parameter $q$. Let us denote by $R_q$ the support of the (homogeneous) multifractal spectrum of $\mu_q$.
Let $h:[0,1]\to R_q$ be a continuous  function. Then for every interval $J\subset [0,1]$, one has
$$\dim\,\{x\in J: h^d_{\mu_q}(x) = h(x) \} = \sup \{d_{\mu_q}(h(x)): x\in J \} .$$
\ET

\sk
We now prove Theorem \ref{theolocalmeas}.

 \sk\sk

 Fix $H>0$, and also $q\in (0,1/2)$.
  If for some $x$ one has \[-  N_0(x) \log_2 p(x)  -  (1-N_0(x)) \log_2 (1-p(x))=h, \]  then there exists 
a real number $h_q(H,x)$ such that  \[ -  N_0(x) \log_2 q  -  (1-N_0(x)) \log_2 (1-q)=h_q(H,x). \]  Since both $p(x)$ and $q$ are strictly less than 
$1/2$, a simple argument entails that the map $h_q(H,x) $ is continuous with respect to $(H,x)$.

Now fix  $x_0\in [0,1]$ and consider the Bernoulli measure with parameter $q=p(x_0)$. Consider the interval $I=B(x_0,r)$. One has
\begin{eqnarray*} 
\{x\in I: h_\mu(x)=H\} &\!\!\!\!  = \!\!\!\!&  \{x\in I: \!-  N_0(x) \log_2 p(x)  -  (1-N_0(x)) \log_2 (1-p(x))=H\}\\
 &\!\!\!\!  = \!\!\!\!& \{x\in I: \!-  N_0(x) \log_2 q  -  (1-N_0(x)) \log_2 (1-q)=h_q(H,x)\}.
\end{eqnarray*}
But this last set has its Hausdorff exactly given by Theorem \ref{th:BQ}, hence
 \begin{eqnarray*} 
\dim \{x\in I: h_\mu(x)=H\} &\!\!\!\!  = \!\!\!\!&   \sup \{d_{\mu_q}( h_q(H,x)): x\in I\}.
\end{eqnarray*}
When $r$ goes to zero, $p(x)$ tends uniformly to $q=p(x_0)$. Hence $h_q(H,x)$ tends to $H$. In particular, the mapping $d_{\mu_q}$ being continuous (real analytic in fact), when $r$ goes to zero one finds that
  \begin{eqnarray*} 
d_\mu(x,H) =  d_{\mu_q}( H) = d_{\mu_{p(x)}}(H).
\end{eqnarray*}

This result can immediately be applied to the case where the mapping $x\mapsto p(x)$ is continuous by part (instead of simply continuous), and can certainly  be adapted when $p$ is c\`adl\`ag.  It would be worth investigating the case where $p$ enjoys less regularity properties.
   
\begin{remark}
Many examples of Cantor set with varying local Hausdorff dimensions have been constructed  \cite{Bae,SchS}, here the key point is that we perform the (global and local) multifractal analysis of measures sitting on these "inhomogeneous" Cantor sets.
\end{remark}

\section{Local  spectrum of stochastic processes}  \label{locspecstoch}

Suppose now that $f$ is a nowhere differentiable function defined on $[0,1]^d$; one can associate with $f$   the  dyadic family $\EEE_f =  \{\Osc_f(3\la)\}_{\la\in \Lambda}$, where the oscillation of $f$ over a set $\omega\subset\Omega$ is
$$\Osc_f(\omega) = \sup\{f(x):x\in \omega\} - \inf\{f(x):x\in \omega\}.$$
Then, it is obvious that  the pointwise \ho exponent \eqref{defpointwise} of $f$ at $x$ is the same as the one defined  by \eqref{exploclimsix} with the dyadic family $\EEE_f$.  
Hence, the previous developments performed in the abstract setting of dyadic functions family holds for non-differentiable functions.

\medskip

We start by giving a simple general probabilistic setting  which naturally leads to
a weak, probabilistic form of homogeneity. 
Let $X$ be a random field on $\RR^d$;  $X$ has stationary increments if 
$ \forall s \in \RR^d, $ the two processes 
\[
x \mapsto Y_s (x) := X(s+x)-X(s) \qquad\mbox{and}\qquad   x \mapsto X(x)
\]
share the same law.  Indeed, this equality in law implies the equality in law of the  linear forms applied to the two processes $Y_s$ and $X$, hence of iterated differences and wavelet coefficients. It  follows that local suprema of iterated differences and of wavelet coefficients computed on dyadic cubes also share the same laws, and Proposition~\ref{regponc} implies that, if $X$ has locally bounded sample paths, then the H\"older exponent has a stationary law. Therefore, the H\"older spectra on dyadic intervals of the same width also share the same law almost surely. As a result, the H\"older spectra  on {\em all} dyadic intervals share the same law.
This leads to the following result.

\BP Let $X$ be a random field on $\RR^d$ with stationary increments. If $X$ has locally bounded sample paths, then
\[
\forall s \quad \as \quad \forall H \qquad d_X (s, H) = d_X (0, H).
\]
\EP


\subsection{Local analysis of the multifractional Brownian motion}\label{sec:mbm}

Let $H$ denote a function defined on $\R^{d}$ with values in a fixed compact subinterval $[a,b]$ of $(0,1)$. We assume that $H$ satisfies locally a uniform H\"older condition of order $\beta\in(b,1)$, that is, $H\in C^{\beta}(\Omega)$ for every open subset $\Omega$ of $\R^{d}$. Now, recall that the multifractional Brownian motion (MBM) with functional parameter $H$ has been introduced in~\cite{Benassi:1997mz,Peltier:1995kx} as the continuous and nowhere differentiable Gaussian random field $B_{H}=\{B_{H}(x),\ x\in\R^{d}\}$ that can be represented as the following stochastic integral
\[
B_{H}(x)=\int_{\R^{d}}\frac{e^{\imath x\cdot\xi}-1}{|\xi|_{2}^{H(x)+d/2}}\,\widehat{dW}(\xi),
\]
where $x\cdot\xi$ denotes the standard inner product, $|\xi|_{2}$ is the usual Euclidean norm, and $\widehat{dW}$ stands for the ``Fourier transform'' of the real-valued white noise $dW$, meaning that for any square-integrable function $f$, one has
\[
\int_{\R^{d}} \widehat{f}(\xi)\,\widehat{dW}(\xi)=\int_{\R^{d}} f(x)\,dW(x).
\]
In particular, the MBM reduces to a fractional Brownian motion when the function $H$ is chosen to be constant. The pointwise regularity of the MBM is well known; as a matter of fact, it has been shown in~\cite{Ayache:2007fk} that
\begin{equation}\label{eq:holdermbm}
\as \quad \forall x\in\R^{d} \qquad h_{B_{H}}(x)=H(x).
\end{equation}
Thus, the H\"older exponent of the MBM is completely prescribed by the function $H$. Our purpose is now to give an illustration to Corollary~\ref{contrex2} above by showing that the multifractal formalism is locally satisfied by almost every sample path of the MBM. To be specific, we shall establish in the remainder of this section the following result which, with the help of Corollary~\ref{contrex2}, enables one to recover~(\ref{eq:holdermbm}).

\begin{prp}\label{prp:MBM}
Let $\EEE_{H}$ denote the dyadic function that is obtained by considering the wavelet leaders of the multifractional Brownian motion $B_{H}$, and assume that the wavelets belong to the Schwartz class. Then, the local scaling function $\tau_{\EEE_{H}}$ satisfies
\[
\as \quad \forall x\in\R^{d} \quad \forall p\in\R \qquad \tau_{\EEE_{H}}(x,p)=H(x)p-d.
\]
\end{prp}

In order to establish Proposition~\ref{prp:MBM}, we shall work with a Lemari\'e-Meyer wavelet basis of $L^{2}(\R^{d})$ formed by the functions $2^{dj/2}\psi^{(i)}(2^{j}x-k)$, see~\cite{Lemarie-Rieusset:1986lr}, and more generally with the biorthogonal systems generated by the fractional integrals of the basis functions $\psi^{(i)}$, namely, the functions $\psi^{(i),h}$ defined by
\[
\widehat{\psi^{(i),h}}(\xi)=\frac{\widehat{\psi^{(i)}}(\xi)}{|\xi|_{2}^{h+d/2}}.
\]
It will also be convenient to consider the Gaussian field $Y=\{Y(x,h),\ (x,h)\in\R^{d}\times(0,1)\}$ given by
\[
Y(x,h)=\int_{\R^{d}}\frac{e^{\imath x\cdot\xi}-1}{|\xi|_{2}^{h+d/2}}\,\widehat{dW}(\xi).
\]
Note, in particular, that $B_{H}(x)=Y(x,H(x))$ for all $x\in\R^{d}$, and that the random field $\{Y(x,h),\ x\in\R^{d}\}$ is merely a fractional Brownian motion with Hurst parameter $h$. By expanding its kernel in the orthonormal basis of $L^{2}(\R^{d})$ formed by the Fourier transforms of the functions $2^{dj/2}\psi^{(i)}(2^{j}x-k)$, and by virtue of the isometry property, the stochastic integral defining $Y(x,h)$ may be rewritten in the form
\[
Y(x,h)=\sum_{i}\sum_{j\in\Z}\sum_{k\in\Z^{d}} \eps^{i}_{j,k} 2^{-hj}\left(\psi^{(i),h}(2^{j}x-k)-\psi^{(i),h}(-k)\right),
\]
where the $\eps^{i}_{j,k}$ form a collection of independent standard Gaussian random variables. It is possible to show that the above series converges uniformly on any compact subset of $\R^{d}\times(0,1)$, see~\cite{Ayache:2005fj}. Moreover, the above decomposition yields the following natural wavelet expansion of the field $B_{H}$:
\begin{equation}\label{eq:wavmbm}
B_{H}(x)=\sum_{i}\sum_{j\in\Z}\sum_{k\in\Z^{d}} \eps^{i}_{j,k} 2^{-H(x) j}\left(\psi^{(i),H(x)}(2^{j}x-k)-\psi^{(i),H(x)}(-k)\right).
\end{equation}
Furthermore, it is shown in~\cite{Ayache:2005fj} that the low-frequency component of $Y$, that is,
\[
\sum_{i}\sum_{j=-\infty}^{-1}\sum_{k\in\Z^{d}} \eps^{i}_{j,k} 2^{-hj}\left(\psi^{(i),h}(2^{j}x-k)-\psi^{(i),h}(-k)\right),
\]
is almost surely a $C^{\infty}$ function in the two variables $x$ and $h$. Hence, the low-frequency component of the MBM, which is obtained by summing only over the negative values of $j$ in~(\ref{eq:wavmbm}), is in $C^{\beta}(\Omega)$ for any open subset $\Omega$ of $\R^{d}$, just as the functional parameter $H$. As $\beta$ is larger than all the values taken by the function $H$, it follows that the pointwise regularity of the MBM is merely given by that of its high-frequency component, that is,
\[
\widetilde{B}_{H}(x)=\sum_{i}\sum_{j=0}^{\infty}\sum_{k\in\Z^{d}} \eps^{i}_{j,k} 2^{-H(x) j}\left(\psi^{(i),H(x)}(2^{j}x-k)-\psi^{(i),H(x)}(-k)\right).
\]
As a consequence, we may consider in what follows the high-frequency component $\widetilde{B}_{H}$ instead of the whole field $B_{H}$. In addition, in view of the regularity of $H$, it follows from standard results on Calder\'on-Zygmund operators (see~\cite{Meyer:1997ve}) and robustness properties of the local scaling functions, $\tau_{\EEE_{H}}$ coincides with the local scaling function of the dyadic family $\widetilde{\EEE}_{H}$ which is obtained by considering the wavelet leaders associated with the wavelet coefficients
\[
c^{i}_{j,k}=\eps^{i}_{j,k} 2^{-H(k2^{-j}) j}.\]
(Recall that in \cite{Jaffard:2004fh}, it is proved that the scaling function is ``robust'', i.e.  does not depend on the smooth enough wavelet basis chosen; furthermore, the arguments of the proof clearly are  local,  so that the local scaling function also is robust.)

Letting $\lambda$ denote the cube corresponding to the indices $i$, $j$ and $k$ as in Section~\ref{obowb}, these coefficients may naturally be rewritten in the form
\[
c_{\lambda}=\eps_{\lambda} 2^{-H(x_{\lambda}) \gen{\lambda}},
\]
where $\eps_{\lambda}$ is the standard Gaussian random variable $\eps^{i}_{j,k}$, $x_{\lambda}$ is the basis point $k2^{-j}$ of the cube $\lambda$ and $\gen{\lambda}$ is its scale $j$. Recall that the wavelet leaders $d_{\lambda}$ are then defined in terms of the wavelet coefficients through~(\ref{eq:defWL}). Finally, for the sake of simplicity and without loss of generality, we shall study the local scaling function $\tau_{\widetilde{\EEE}_{H}}$ only on the open set $(0,1)^{d}$, so that we only have to consider the dyadic subcubes of $[0,1)^{d}$.

Let us now establish a crucial lemma concerning the behavior on the subcubes of $[0,1)^{d}$ of the new dyadic family $\widetilde{\EEE}_{H}$.

\BL\label{lem:WLmbm}
With probability one, for any dyadic cube $\lambda\subset [0,1)^{d}$ with scale $\gen{\lambda}$ large enough,
\[
\frac{1}{\gen{\lambda}^{3H(x_{\lambda})}}\leq 2^{H(x_{\lambda}) \gen{\lambda}} d_{\lambda}\leq 2\gen{\lambda}.
\]
\EL

\begin{proof}
We begin by the proving the lower bound. For any proper dyadic subcube $\lambda$ of $[0,1)^{d}$ with scale $\gen{\lambda}=j$, we have
\[
\prob(d_{\lambda}\leq\gen{\lambda}^{-3H(x_{\lambda})} 2^{-H(x_{\lambda}) \gen{\lambda}})
=\prod_{\lambda'\subset 3\lambda}\prob(|\eps_{\lambda'}|\leq\gen{\lambda}^{-3H(x_{\lambda})} 2^{H(x_{\lambda'}) \gen{\lambda'}-H(x_{\lambda}) \gen{\lambda}}).
\]
Let $l(j)=j+\lceil (2/d)\log_{2}j \rceil$, where $\lceil\,\cdot\,\rceil$ denotes the ceiling function and $\log_{2}$ the base two logarithm. Considering in the above product only the subcubes $\lambda'\subset 3\lambda$ with scale $\gen{\lambda'}$ equal to $l(j)$, and using the elementary fact that the modulus of a standard Gaussian random variable is bounded above by $t$ with probability at most $t$, we deduce that
\[
\prob(d_{\lambda}\leq\gen{\lambda}^{-3H(x_{\lambda})} 2^{-H(x_{\lambda}) \gen{\lambda}})
\leq\prod_{\lambda'\subset 3\lambda \atop \gen{\lambda'}=l(j)} \gen{\lambda}^{-3H(x_{\lambda})} 2^{H(x_{\lambda'}) \gen{\lambda'}-H(x_{\lambda}) \gen{\lambda}}.
\]
Moreover, the function $H$ satisfies locally a uniform H\"older condition of order $\beta$, so there exists a real $C>0$ that does not depend on $\lambda$ such that
\begin{equation}\label{eq:holdHlambda}
\forall \lambda'\subset 3\lambda \qquad |H(x_{\lambda'})-H(x_{\lambda})|\leq C 2^{-\beta j}.
\end{equation}
Combined with the observation that there are at least $j^{2}$ subcubes $\lambda'\subset 3\lambda$ such that $\gen{\lambda'}=l(j)$, this implies that
\[
\prob(d_{\lambda}\leq\gen{\lambda}^{-3H(x_{\lambda})} 2^{-H(x_{\lambda}) \gen{\lambda}})
\leq \left(j^{-3H(x_{\lambda})} 2^{H(x_{\lambda})(l(j)-j)+C l(j) 2^{-\beta j}}\right)^{j^{2}}.
\]
Given that the function $H$ is valued in the interval $[a,b]$, we infer that
\[
\prob(d_{\lambda}\leq\gen{\lambda}^{-3H(x_{\lambda})} 2^{-H(x_{\lambda}) \gen{\lambda}})
\leq \left(j^{(2/d-3)a} 2^{b+C l(j) 2^{-\beta j}}\right)^{j^{2}}.
\]
The right-hand side is clearly bounded above by $e^{-j^{2}}$ when $j$ is larger than some integer $j_{0}$, so that
\[
\sum_{\lambda\subset [0,1)^{d} \atop \gen{\Lambda}\geq j_{0}} \prob(d_{\lambda}\leq\gen{\lambda}^{-3H(x_{\lambda})} 2^{-H(x_{\lambda}) \gen{\lambda}})
\leq\sum_{j\geq j_{0}} 2^{dj}e^{-j^{2}}<\infty,
\]
and we deduce the required lower bound from the Borel-Cantelli lemma.

In order to establish the upper bound, let us begin by observing that with probability one, for any dyadic cube $\lambda\subset [0,1)^{d}$ with scale $\gen{\lambda}=j$ large enough, $|\eps_{\lambda}|\leq j$. This follows again from the Borel-Cantelli lemma, together with the fact that
\[
\prob(|\eps_{\lambda}|>j)=2(1-\Phi(j))\leq\frac{e^{-j^{2}/2}}{j}\sqrt{\frac{2}{\pi}},
\]
which itself follows from standard estimates on the asymptotic behavior of the cumulative distribution function $\Phi$ of the standard Gaussian distribution. Now, along with~(\ref{eq:holdHlambda}), this implies that for $\gen{\lambda}=j$ large enough,
\[
d_{\lambda}\leq \sup_{\lambda'\subset 3\lambda} \gen{\lambda'} 2^{-(H(x_{\lambda})-C 2^{-\beta j}) \gen{\lambda'}}
=j2^{-(H(x_{\lambda})-C 2^{-\beta j}) j}
\leq 2j2^{-H(x_{\lambda})j},
\]
and the required upper bound follows.
\end{proof}

We may now finish the proof of Proposition~\ref{prp:MBM}. To this end, let $x\in (0,1)^{d}$ and $r>0$ such that $\Omega=B(x,r)\subset (0,1)^{d}$. Then, owing to Lemma~\ref{lem:WLmbm}, the structure function of the dyadic function $\widetilde{\EEE}_{H}$ on $\Omega$, which is defined by~(\ref{defstrucfunc}), satisfies
\begin{equation}\label{eq:encstrucmbm}
\sum_{ \la \in \La_j^\Ome} \left(\frac{2^{-H(x_{\lambda})j}}{j^{3H(x_{\lambda})}}\right)^{p} \leq S_j (\Ome, p)\leq \sum_{ \la \in \La_j^\Ome} \left(2j2^{-H(x_{\lambda})j}\right)^{p}
\end{equation}
for $j$ large enough and $p\geq 0$. Given that $H$ satisfies locally a uniform H\"older condition of order $\beta$, there exists a real $C>0$ that depends on neither $x$ nor $r$ such that $|H(x_{\lambda})-H(x)|\leq C r^{\beta}$ for all dyadic cubes $\la\subset\Ome$. In addition, the cardinality of $\La_j^\Ome$ is comparable with $r^{d}2^{dj}$. Thus, there is a constant $C'>0$ such that
\[
\frac{r^{d}2^{dj}}{C'} \left(\frac{2^{-(H(x)+Cr^{\beta})j}}{j^{3(H(x)+Cr^{\beta})}}\right)^{p} \leq S_j (\Ome, p)\leq C' r^{d}2^{dj} \left(2j2^{-(H(x)-Cr^{\beta})j}\right)^{p}.
\]
It follows that the scaling function of $\widetilde{\EEE}_{H}$ on $\Omega$ satisfies
\[
(H(x)-Cr^{\beta})p-d\leq\tau_{\widetilde{\EEE}_{H}}^\Ome (p)\leq (H(x)+Cr^{\beta})p-d.
\]
Letting $r$ go to zero, we may finally conclude that $\tau_{\EEE_{H}}(x,p)=H(x)p-d$ for all $p\geq 0$ and $x\in(0,1)^{d}$. The same approach still holds for the negative values of $p$ except that the inequalities have to be reversed in~(\ref{eq:encstrucmbm}) and in the subsequent estimates as well. Proposition~\ref{prp:MBM} follows.

 \subsection{ A Markov process with a varying local multifractal spectrum}\label{sec:markov}

  In this section we reinterpret the results of  \cite{Barral:2010vn} in terms of local spectrum.
A quite general class of one-dimensional Markov processes consists of
stochastic differential equations (S.D.E.) with jumps. Recall that such a process is the sum of a Brownian motion and a pure jump process.  
We will assume in the following that the process has no Brownian part;  indeed, since  Brownian motion is mono-H\"older, its consequence on the spectrum is straightforward to handle: it eliminates H\"older exponents  larger than $1/2$  and, eventually adds a point at $(1/2, 1)$. Thus the Markov processes that will be studied are jumping S.D.E. without Brownian and drift part, starting e.g. from $0$,
and with jump measure $\nu(y,du)$ (meaning that 
when located at $y$, the process 
jumps to $y+u$ at rate $\nu(y,du)$). Again, since this is a "toy" model, we will make additional simplifying assumptions: Namely that   the process  is increasing 
(that is, $\nu(y,(-\infty,0))=0$ for all $y \in \R$). 
Classically, a necessary condition for the process to be well-defined is that
$\int_0^\infty u \,\nu(y,du)<\infty$.

\sk

If ${\nu }$ is chosen so that the index $ \beta_{{\nu }(y,.)}$ is constant with respect to $y$,  then one expects that  the local multifractal spectrum $ d_M(t,h)$  of the process $M=(M_t)_{t\geq 0}$ will be deterministic and independent of $t$.
Hence,  the index of the jump measure
will depend on the  value
$y$ of the process. The most natural example of such a situation 
consists in choosing
\begin{equation*}
{\nu_\gamma}(y,du):= \gamma(y) u^{-1-\gamma(y)} \indiq_{[0,1]}(u)du,
\end{equation*}
for some function $\gamma: \R \mapsto (0,1)$.  The lower exponent of this family of measures is 
$$
\forall \, y\geq 0, \ \ \\ \ \ \  \beta_{{\nu_\gamma}(y,.)}=\gamma(y).
$$

In \cite{Barral:2010vn}, the following assumption is made

\begin{equation*}
({ \mathcal H} ) \hspace{5mm} \begin{cases}  \mbox{There exists $\e>0$  such that $\gamma:[0,\infty)  \longmapsto [\ep,1-\ep]$ }
\\\mbox{ is a Lipschitz-continuous strictly increasing  function.}\end{cases}
\end{equation*}

It is relatively clear that the assumptions can be relaxed, and that many classes of Markov processes could be further studied. An interesting subject to investigate is  the range of functions $\gamma$ that could be used in the construction. For a process, $M=(M_t)_{t\geq 0}$, one sets $\Delta M_t=M_t - M_{t-}$, 
where $M_{t-} = \displaystyle\lim_{s\to t, \, s<t} M_s$

\begin{prp}\label{exi}\cite{Barral:2010vn}
Assume  that $({ \mathcal H} )$ holds.  There exists a strong Markov process
$M=(M_t)_{t\geq 0}$   starting from $0$,   increasing and c\`adl\`ag (i.e. right-continuous, with a left limit),
and with generator $\cL$
defined for all $y \in [0,\infty)$ and for any function 
$\phi: [0,\infty) \mapsto \R$ 
Lipschitz-continuous by
\begin{equation}\label{giok}
\cL \phi(y) = \int_0^1 [\phi(y+u)-\phi(y)] 
{\nu_\gamma}(y,du).
\end{equation}
Almost surely, this process is  continuous except on  a countable number of jump times. Denote
by $\cJ $ the set of its jump times, that is $\cJ=\{t >0: \Delta M(t)\ne 0\}$.
Finally, $\cJ $ is dense in $[0,\infty)$.
\end{prp}
  
This representation of $M$ is useful for its local regularity 
analysis.

%


%

\smallskip

The following  theorem  of \cite{Barral:2010vn} summarizes the multifractal  features of $M$.

\begin{thm}\label{mr1}
Assume $(\cH)$ and consider the process $M$ 
constructed  in Proposition \ref{exi}. 
Then, the following properties hold  almost surely.  
\begin{enumerate}
\item[(i)]
For every $t\in (0,\infty)\backslash \cJ$, 
the local spectrum of $M$ at $t$  is  given by 
\begin{equation}\label{localspecM}
d_M(t,h)  = \begin{cases} h \cdot  \gamma(M_t)     &  \mbox{ if } \  0\leq h\leq 1/\gamma(M_t)  , \\  
- \infty &  \mbox{ if } \ h > 1/\gamma(M_t),\end{cases}
\end{equation}
while for $t\in \cJ$, 
\begin{equation}\label{localspecM2}
d_M(t,h)  = \begin{cases} h \cdot \gamma(M_t)      &  \mbox{ if } \  0\leq h<1/\gamma(M_t)  , \\ 
h \cdot  \gamma(M_{t-})    &  \mbox{ if } \  h\in [1/\gamma(M_t), 1/\gamma(M_{t-})]  ,\\ - \infty &  \mbox{ if } \ h > 1/\gamma(M_{t-}).\end{cases}
\end{equation}
\item[(ii)]
The spectrum of $M$ on any interval $I=(a,b)\subset (0,+\infty)$  is  
\begin{eqnarray}
\label{mfs1} \forall  h \geq 0, \ \ \  d_M(h) & =&          \sup\Big \{h\cdot  \gamma(M_t): \, t\in I, \ h\cdot  \gamma(M_t) <1\Big\} \\
\label{mfs1.5} & =&  \sup\Big \{h\cdot  \gamma(M_\sm): \, s \in \cJ\cap I, \ h\cdot  \gamma(M_\sm) <1\Big\} .\end{eqnarray}
In \eqref{mfs1} and \eqref{mfs1.5}, we adopt the convention that $\sup \emptyset =-\infty$.
\end{enumerate}
\end{thm}

\begin{figure}[h]
\includegraphics[width=6.cm,height = 4cm]{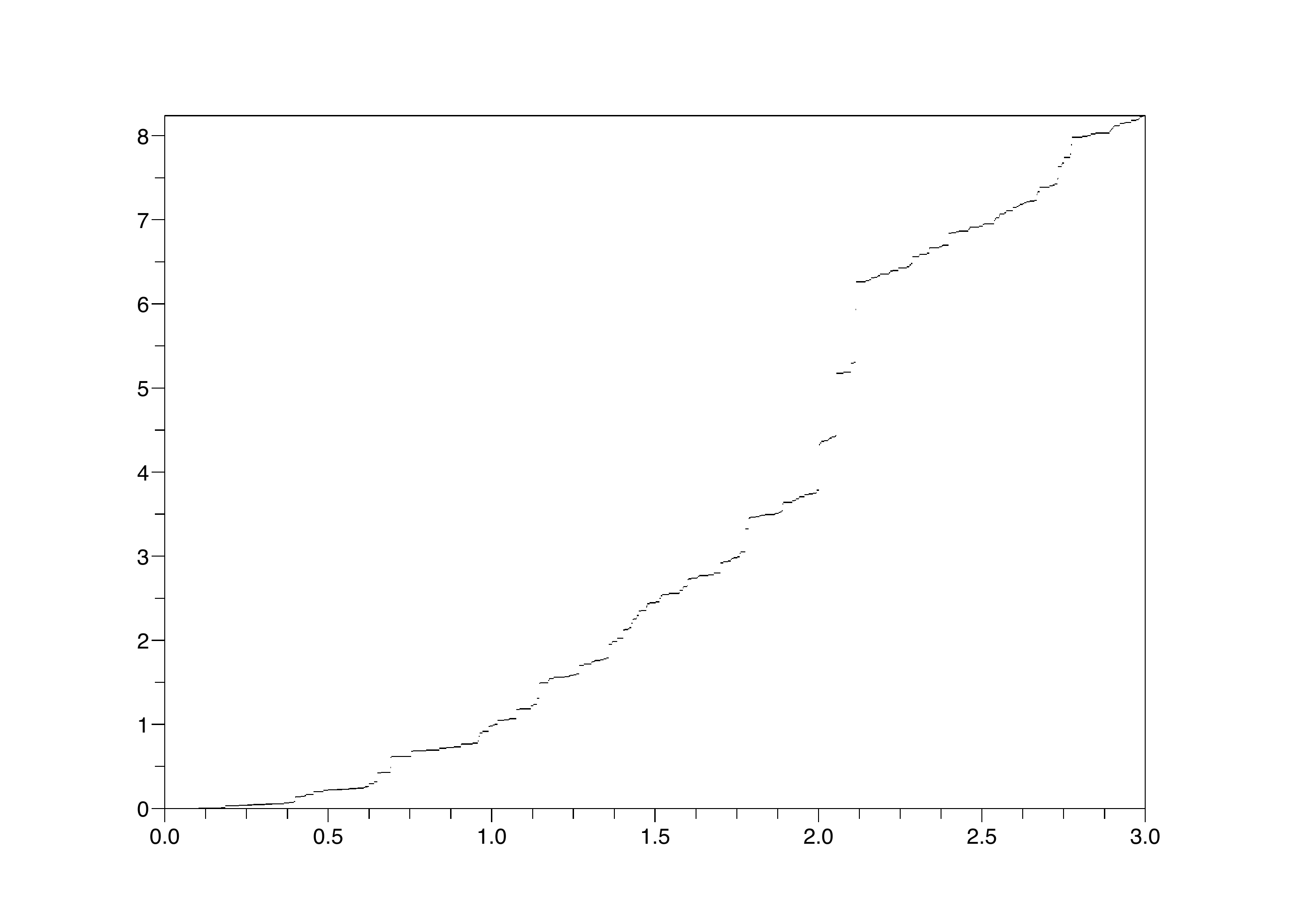}  
\includegraphics[width=6.cm,height = 4cm]{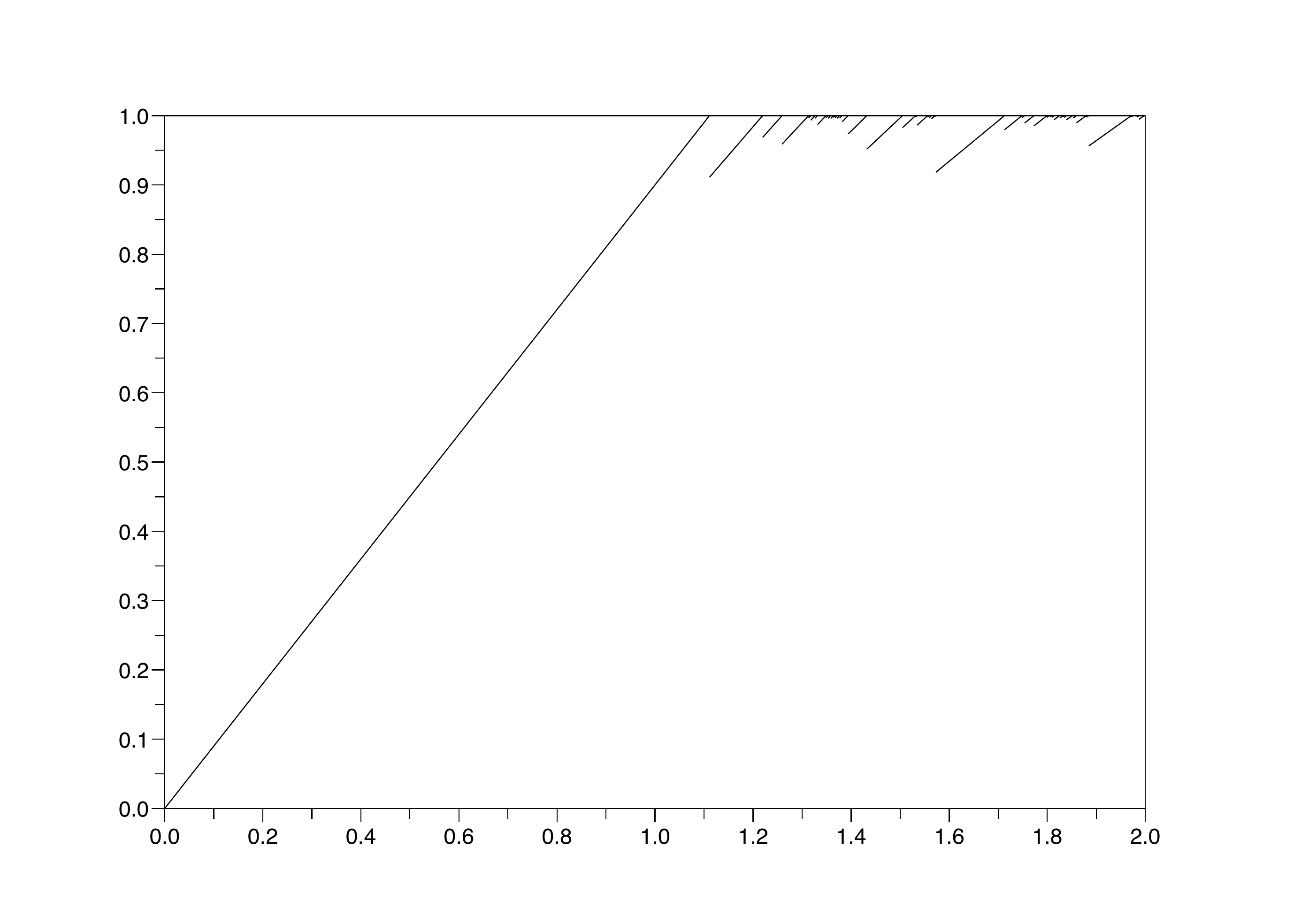}  
\includegraphics[width=6.cm,height = 4cm]{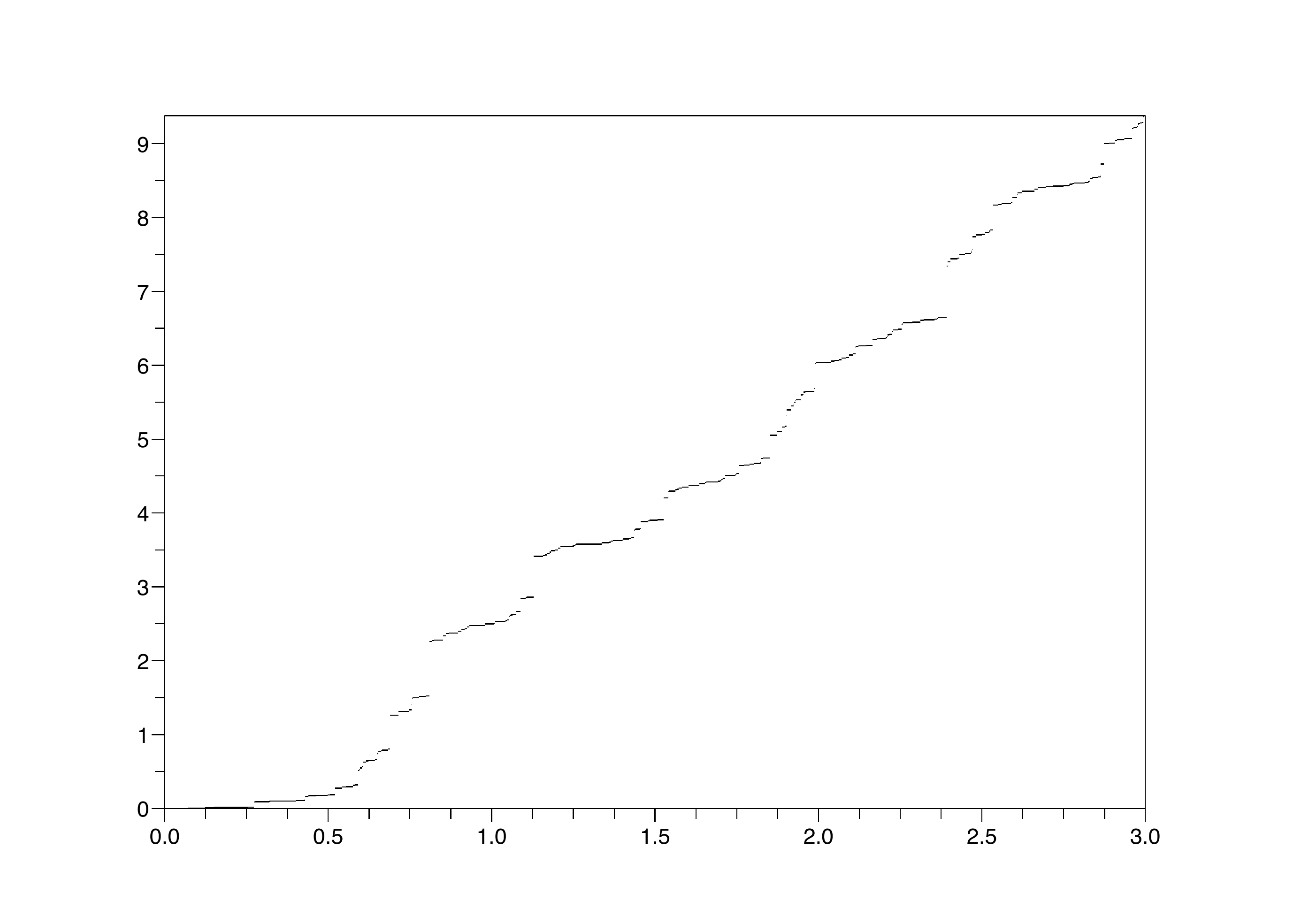}  
\includegraphics[width=6.cm,height = 4cm]{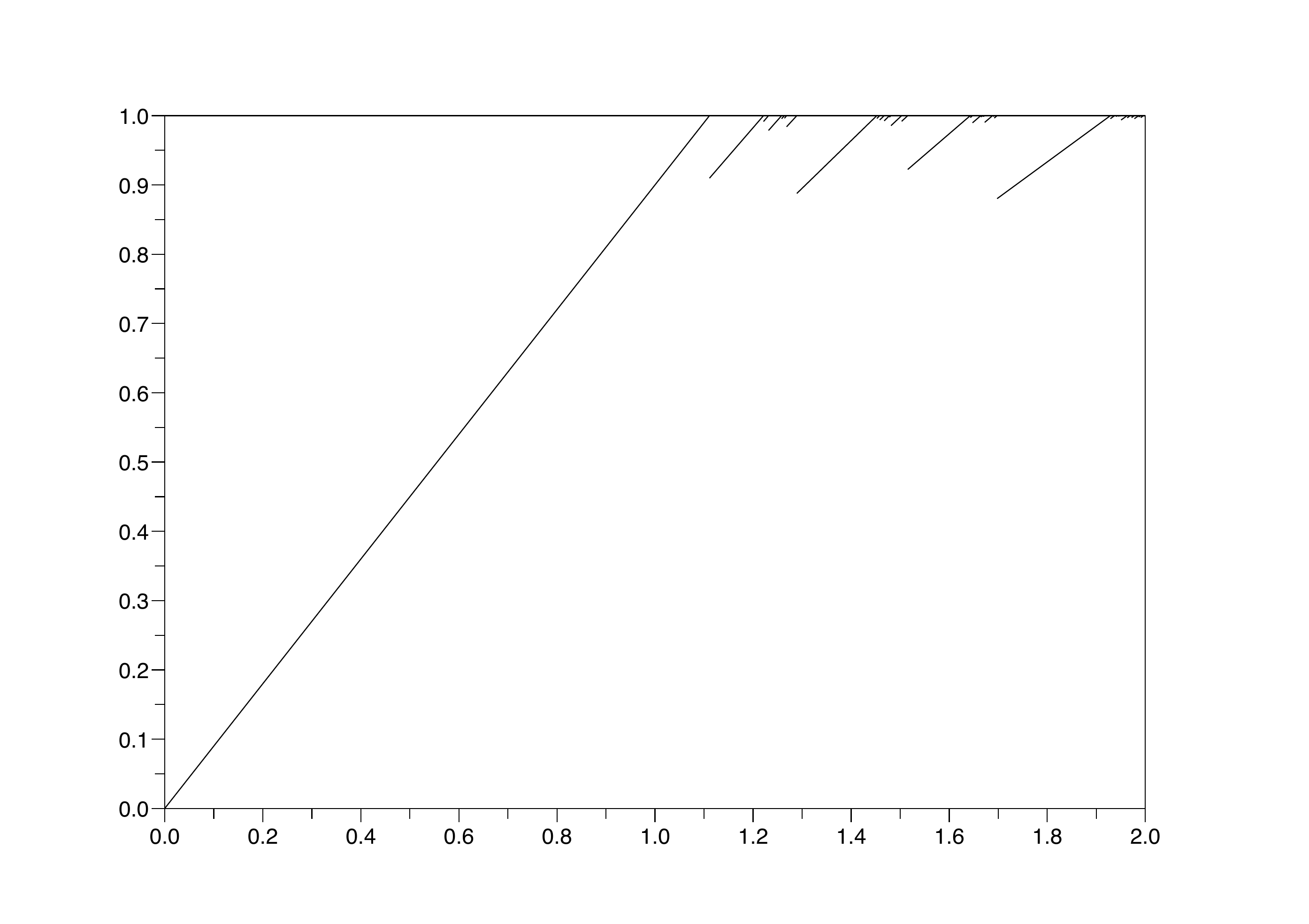} 
\caption{Two sample paths of the stochastic process $M$ built using the function $\gamma(y):=\min(1/2+y/4,0.9)$. On the right hand-side are plotted the 
theoretical spectra $D_M([0,3],.)$.}
\end{figure}

As can be seen from the definition of the local multifractal spectrum, in order to prove Theorem \ref{mr1}, it is enough to show \eqref{mfs1.5}. Indeed, \eqref{localspecM} simply follows from considering the limit of \eqref{mfs1.5} when the interval $I$ is the centered ball $B(t,r)$ and letting $r$ tend to zero. 

  \sk\sk

Formula \eqref{mfs1.5} is better understood when plotted: 
for every $s\in I \cap \cJ$, plot a segment whose endpoints are $(0,0)$ and  
$(1/\gamma(M_{s-}),1)$ (open on the right), and take the supremum to get 
$D_M(I,.)$. Sample paths of the process $M$ and their associated spectra are given in Figure 1.


The formulae giving the local and global spectra are based on the computation of the pointwise \ho exponents at all times $t$. The value of the pointwise \ho exponent of $M$ at $t$ depends on two parameters: the value of the process $M$ in the neighborhood of $t$, and the approximation rate of $t$ by the set of jumps $\cJ$. In particular, the following properties holds a.s.,
\begin{eqnarray*}
\mbox{for   every $t\geq 0$, } \ &&h_M(t) \leq 1/\gamma(M_t),\\
\mbox{for Lebesgue-almost every $t$, }&& h_M(t)= 1/\gamma(M_t),\\
\mbox{for every $\kappa\in (0,1)$, }&& 
\dim_H \{t\geq 0: \; h_M(t)= \kappa/\gamma(M_t) \}=\kappa.
\end{eqnarray*}

The relevance of the local spectrum in this context is thus obvious: depending on the local value of $M$, the pointwise \ho exponents change, and so is the (local) multifractal spectrum.


It is worth emphasizing that, as expected from the construction of the process $M$, the local spectrum \eqref{localspecM} at any point $t>0$ essentially
coincides with that of a stable L\'evy subordinator of index $\gamma(M_t)$. This local comparison is strengthened by the following theorem, which  proves the existence of tangent processes for $M$ (which are L\'evy stable subordinators).

\begin{prp}\label{tangent}
We denote by 
$\mathcal{F}_t:=\sigma(\{N(A), A\in \mathcal{B}([0,t] \times [0,\infty))\})$. Let $t_0\ge 0$ be fixed. Conditionally on $\mathcal{F}_{t_0}$, 
the family of processes 
$\displaystyle \Big (\frac{M_{t_0+\alpha t}-M_{t_0}}
{\alpha^{1/\gamma(M_{t_0})}}\Big )_{t\in [0,1]}$ converges in law, 
as $\alpha\to 0^+$, to a stable L\'evy 
subordinator with L\'evy measure $\gamma(M_{t_0}) u^{-1-\gamma(M_{t_0})} du$.
Here the Skorokhod space of c\`adl\`ag functions on $[0,1]$ is endowed
with the uniform convergence topology.
\end{prp}

 Observe that for all $s\in \cJ$, all $h\in (1/\gamma(M_s),1/\gamma(M_\sm)]$,  
$d_M(h) = h \cdot \gamma(M_\sm)  $. Thus the 
spectrum $d_M$ of $M$ on an interval $I$ is a straight line on all segments of the form 
$(1/\gamma(M_s),1/\gamma(M_\sm)]$, $s\in \cJ \cap I$. By the way, this spectrum, 
when viewed as a map from $\R_+$ to $\R_+$,  is very irregular,
and certainly multifractal itself. 
This is in sharp contrast with  the spectra usually obtained, which are most of the time concave or (piecewise) real-analytic. 
Hence, the difference between the global and the local multifractal spectra is stunning: While $d_M$ is very irregular, $d_M(t,\cdot)$ is a straight line.

 This example naturally leads to  the following  open problem, which would express that a  natural compatibility holds for  local multifractal analysis: Find general conditions under which  a stochastic process $X$ which has a tangent process at  a point $x_0$ satisfies that  the  multifractal spectrum of the tangent process  coincides with the local spectrum of $X$ at $x_0$.

\section{Other regularity exponents characterized by dyadic families}  \label{other}

Other exponents than those already mentioned fit in the general framework given by Definition \ref{defexpabstra} and therefore the  results supplied by multifractal analysis can be applied to them.  We now list a few of them.



Pointwise H\"older regularity  is pertinent only if applied to locally bounded functions.  An extension of  pointwise  regularity fitted to functions that are  only assumed to  belong to  $L^p_{loc}$ is sometimes required: The corresponding notion  was introduced by  Calder\'on and Zygmund in  1961, see
\cite{Calderon:1961kx},  in order  to obtain pointwise regularity results for elliptic PDEs.

\BD \label{def2}  Let   $p \in [1, +\infty)$ and $\al > -d/p$.  Let $f\in L^p_{loc} (\Ome) $, and $x_0 \in \Ome$;  $f $  belongs to 
   $T^p_\al (x_0)$ if   there exist $\ 
C >0$  and a polynomial $P$ of degree less than   $\al$  such that, for $r$ small enough,   
\begin{equation}  \label{deux}  \left( \frac{1}{r^d} \int_{B(x_0,r)} | f(x) -P( x-x_0)|^p dx 
\right)^{1/p}
\leq C r^\al.
\end{equation}
The
$p$-exponent of  $f$ at $x_0$ is \[  h^{p}_f
(x_0)= \sup \{ \al : f \in T^p_\al (x_0)\}.  \]
\ED

{ \bf Remarks:} 

\begin{itemize}
\item The normalization chosen in (\ref{deux}) is such that cusps $|x-x_0|^\al$ (when $\al \notin 2\N$) have an H\"older and a $p$-exponent which take the same value $\alpha$ at $x_0$. 
\item  The H\"older exponent corresponds to the case $p =+\infty$.  
\item We only define lower exponents here: Upper exponents could also be defined in his context, by considering local $L^p$ norms of iterated differences. 
\item Definition \ref{def2}   is a natural substitute for pointwise H\"older regularity when
 functions in 
$L^p_{loc}$ are considered. In particular, the $p$-exponent can take negative values  down to $-d/p$, and typically  allows to take into account behaviors which are locally of the form 
\begin{equation}\label{eq:defHecke}
\frac{1}{|x-x_0|^\gamma }\;\;   \mbox{ for}\;\; \gamma < d/p, 
\end{equation}
\end{itemize}

A pointwise regularity exponent  associated with tempered distributions has been introduced by Y. Meyer: The { \em weak scaling exponent} (see \cite{Mey1}, and also \cite{AJL} for a multifractal formalism based on this exponent).  It can also be  interpreted as  a limit case of other exponents for distributions, which can be related with the  H\"older exponent; let us briefly recall how this can be done. 

Let $f$ be a tempered distribution defined over $\R^d$.  One can define fractional primitives of order $s$ of $f$ in the Fourier domain by 
\[ \widehat{f^{(-s)}} (\xi ) = (1+ | \xi |^2)^{s/2} \hat{f} (\xi) .\]
 Since  $f$ is of finite order, for $s$ large enough, $f^{(-s)}$ locally belongs  to $L^p$ (or $L^\infty$). It follows that one can define regularity exponents of distributions  through  $p$-exponents (or H\"older exponents) of a	 fractional primitives of large enough order. If $f$ is only defined on a domain $\Ome$,  one can still define the same exponents  at $x_0\in \Ome$ by  using a function $g \in {\mathcal D} (\R^d)$ such that  $g$ is supported inside $\Ome$ and $ g(x) =1$ in a neighborhood of $x_0$; then $fg$ is  a tempered distribution defined on $\R^d$ and the  exponents  of $(fg)^{(-s)}$ at  $x_0$ clearly do not depend on the choice of $g$.

 Let $f$ be a tempered distribution defined on a open domain. 
 Denote by $h^s_f (x)$ the H\"older exponent of $f^{(-s)}$ (which is thus canonically well defined for $s$ large enough. By definition, the weak scaling exponent of $f$ at $x$ is 
 \[ { \mathcal W}_f (x) = \lim_{s \rightarrow + \infty}  \left( h^s_f (x) - s\right) \]
 (note that  the limit always exists because the quantity considered is an increasing function of $s$). 
 We will not deal directly with this exponent because it does not directly fit in the framework given by  Definition \ref{defexpabstra}. But we will rather consider the following intermediate framework.

 \BD Let $f$ be a tempered distribution defined on a  non-empty open set $\Ome \subset \R^d$. Let $ p \geq 1$ and  $s$ be large enough so that $f^{(-s)}$ belongs to $L^p$ in  a neighborhood of $x_0$. The fractional $p$-exponent of order $s$ of $f$ at $x_0$ 
 is defined by 
 \[Ê  h^{p,s}_f (x_0) = h^{p}_{f^{(-s)}}(x_0) \]
 (using the convention  $h^{\infty}_f = h_f $).
 \ED
 
 Note that, in practice, the standard way to perform  the multifractal analysis of data that are not locally bounded  is to deal with the exponent  $h^{\infty,s}_f $, where $s$ is chosen large enough so that $f^{(-s)} \in L^\infty_{loc}$, i.e. it consists in first performing a fractional integration, and then  a standard multifractal analysis based on the H\"older exponent, see \cite{AJW1} and references therein.

Similarly,  in the function case, if the pointwise regularity exponents are small enough, they can be recovered for the oscillation of $f$. Recall that the  oscillation of $f$ of order $l$ on a convex set $A$ is defined through  conditions on the finite differences of the function $f$, denoted by  $\Delta^M_h f$:   The first order difference of $f$ is  \[ (\Delta_h ^1 f)(x) = f(x+h) -f(x). \] If $l >1 $, the differences of order  $l$ are defined recursively by \[ (\Delta^l_h f)(x) = (\Delta^{l-1}_h f)(x+h)-(\Delta^{l-1}_h f)(x).\]
Then
\[ \Osc^l_f (A) = \sup_{x, x+lh \in A} \left| (\Delta^l_h f)(x) \right| .\]
One easily checks that the H\"older exponent can be derived for the oscillation  on the cubes $3 \la$. Let  $f$  be  locally bounded    on an open set $\Ome$.

 If $l  >  h_f(x_0) $, then 
\BE \label{exploclimbis1} \forall x_0 \in \Ome  , \hspace{1cm} h_f(x_0) = \liminf_{j \rightarrow + \infty}   \frac{ \log 
  \Osc^l_f (  3\la_{j}(x_0) ) }{\log  \, 2^{-j}  }  .
\end{equation}

Recall also Proposition \ref {regponc} which    allows to derive numerically  the H\"older exponent  by a log-log plot regression  bearing on the  the  $ d_{ \la_j (x_0)}$  when $j \rightarrow +
\infty$, see  \cite{Jaffard:2004fh}.

However, in contradistinction with the measure case,  a similar formula does not hold for the upper H\"older exponent, see \cite{Claus} where partial results in this direction and counterexamples are worked out.  \\

We now turn to the wavelet characterization of  the $p$-exponent. 
 We will assume that  $f$ locally belongs  to $L^p$, with slow $L^p$-increase, i.e. satisfies 
\[ \exists C, N >0  \hspace{1cm} \int_{ \Ome \cap B(0, R) } |  f(x) |^p dx  \leq C (1+ | R|)^N . \] 

  In the following, when dealing with  the $ T^p_\al $  regularity of a  function $f$, we will always assume that, if $f$ is  defined on an unbounded set $\Ome$, then it has  slow $L^p$-increase, and, if $\Ome \neq \RR^d$, then  the wavelet basis used is compactly supported.

\BD Let $f \in L^p_{loc} (\Ome)$, and 
let $\psi_\la$ be a given wavelet basis. 
 The {  local square function  }  of $f$ is
\[ S_{f ,\la } (x)= \left( \sum_{ \la' \subset 3 \la} | \clap |^2 1_{\la'} (x)  \right)^{1/2}, \]
and the $p$-leaders are defined by $  d^p_\la =  2^{dj/p} \parallel  S_{f ,\la }  \parallel_p . $ 
\ED

The following result of \cite{JaCiel}  yields a wavelet characterization of  the $p$-exponent which is similar to   (\ref{leaderregpon}). 

\BP \label{prop:caracpexp} Let $p \in (1,\infty)$ and   $f\in L^p$.     Assume that the wavelet basis used is $r$-smooth with 
$r >  h^p_f(x_0) +1$. Then
\BE \label{exploclimcinq} h^p_f(x_0) =  \liminf_{j \rightarrow + \infty}   \frac{ \log 
    d^p_{\la_j (x_0)}   }{\log  \, 2^{-j}  }  .
\end{equation}
\EP

Recall that the ``almost-diagonalization'' principle for  fractional integrals on wavelet bases states that,  as regards H\"older regularity, function spaces or scaling functions,  one can consider that a fractional integration just acts as  if it were diagonal on a wavelet basis,  with coefficients $2^{-sj}$ on $\psi_\la$.  This  rule of thumb  is justified by the fact that  a fractional integration actually is the product of such a diagonal operator and of an invertible Calderon-Zygmund operator $A$ such that $A$ and $A^{-1}$ both belong to the { \em  Lemari\'e algebras}  ${ \mathcal M}^\gamma$, for a $\gamma$ arbitrarily large (and which depends only on the smoothness of the wavelet basis) see \cite{Mey, Meyer:1997ve} for the definition of the Lemari\'e algebras and for the result concerning  function spaces and \cite{Jaffard:2004fh} and references therein for H\"older regularity, function spaces or scaling functions.  
   
  It follows from Proposition \ref{prop:caracpexp} , and the ``almost-diagonalization'' principle for  fractional integrals on wavelet bases, that the exponent $h^{p,s}_f (x_0) $ can be obtained  as follows.

\BC Let $p \in (1,\infty)$ and   $f\in L^p$.  Let 
\[ S^s_{f ,\la } (x)= \left( \sum_{ \la' \subset 3 \la} | 2^{-sj'} \clap |^2 1_{\la'} (x)  \right)^{1/2} \;\; \mbox{ and} \;\; \;\;   d^{p,s}_{\la} 
=  2^{dj/p} \parallel  S^s_{f ,\la }  \parallel_p .  \]
 Then, if the wavelet basis is  $r$-smooth with 
$r >  h^p_f(x_0) + s + 1$, then 
\BE \label{exploclimsixt} h^{p,s}_f (x_0)  =  \liminf_{j \rightarrow + \infty}   \frac{ \log 
   d^{p,s}_{\la_j (x_0)}  }{\log  \, 2^{-j}  } .
\end{equation}
\EC


\section{A  functional analysis point of view}   \label{functional} 

 
  
  \subsection{Function space  interpretation: Constant regularity  }

 If $p >0$, the scaling function has  a function space interpretation, in terms of { \em discrete Besov spaces} which we now define.  Recall that the elements of a dyadic family are always non-negative.

 \BD \label{defvin} Let $s \in \RR$ and $p \in \RR $. A dyadic function     $\EEE$ belongs to  $\bsp(\Ome )$   if 
\BE  \label{dos}  \exists C, \;  \forall j ,  \hspace{1cm}
2^{-dj} {\sum_{\la \in \La^\Ome_j} }^{ \ast  }(e_\la )^p \leq C\cdot 2^{-spj}.\EE
If $p = + \infty$, a dyadic function   $\EEE$ belongs to  ${\bf b}^{s, \infty}_\infty (\Ome ) $  if 
\BE  \label{dos2}  \exists C \;  \forall \la  : \hspace{4mm}
 e_\la   \leq C\cdot 2^{-sj}.\EE
\ED

Note that, if $p >0$, this condition (if applied to the moduli of the coefficients) defines a vector space. It is a Banach space if $p \geq 1$, and a quasi-Banach space if $0 < p <1$; recall that, in a  quasi-Banch space, the triangular inequality is replaced by the weaker condition : 
\[ \exists C, \; \forall x, y, \hspace{1cm} \parallel x+y \parallel \leq C ( \parallel x \parallel + \parallel y \parallel) . \]  Definition  \ref{defvin}  yields  a  function space interpretation  to  the scaling function when $p >0$.  It is classical in this context to rather consider the scaling  function
  $$\eta_\EEE(p)= \tau_\EEE(p) -d.$$
Then,  if $\Ome$ is a bounded set,  
\[ \forall  p \in \RR , \hspace{1cm}  \eta_{\EEE}^\Ome (p)  =   \sup \{ s: \EEE \in   {\bf b}^{s/p, \infty}_p (\Ome )  \} ; \] 
and,  if $\Ome$ is unbounded, then the function space interpretation  is the same, using the precaution supplied by (\ref{omenotboun}). Additionally, 
\[ \Hmin = \sup \left\{ s: \EEE \in  {\bf b}^{s, \infty}_\infty (\Ome )  \right\}  .  \]

The terminology of ``discrete 
  Besov  spaces''  is justified by the fact that, if the $e_\la$ are  wavelet coefficients, then  (\ref{dos})  and (\ref{dos2}) are  the wavelet characterization of the ``classical'' Besov spaces $B^{s, \infty}_p (\RR^d) $ of functions (or distributions) defined on $\RR^d$;
   therefore each wavelet decomposition establishes an isomorphism between  the space $\bsp (\RR^d)  $ and the space  $B^{s, \infty}_p (\RR^d) $, see \cite{Mey}. Note that, when $p = \infty$, these Besov spaces  coincide with the H\"older spaces $C^s (\RR^d)  $, so that, when the $(e_\la)$ are wavelet coefficients, then  the uniform regularity exponent has the following interpretation
   \[  \Hmin =  \sup \{ s: \EEE \in   C^s  (\RR^d)  \} . \] 
 
   
  In the measure case,  H. Triebel showed that the discrete Besov conditions bearing on the $\mu (3 \la)$ can also be related with  the Besov regularity of the measure $ \mu$, see \cite{Triebel}: 
  \[ \mbox{ If } s <d, \hspace{8mm}  (\mu (3 \la )) \in \bsp   (\RR^d)  \Longleftrightarrow  \mu \in  B^{s-d, \infty}_p (\RR^d) .  \]
  In the case where $p = + \infty$,  uniform regularity  gives an important information concerning the sets  $A$ such that $\mu (A) >0$, as a consequence of the   { \em mass distribution principle}, see Section \ref{sec:locmulanal}:
 Since this estimate precisely means  that  the sequence $(e_\la ) = (\mu (3 \la)) $ belongs to $ {\bf b}^{s, \infty}_\infty (\Ome ) $, it follows that, if $A\subset \Ome$  and if a measure $\mu$   satisfies $\mu (A) >0$, then 
 \BE \label{regunifmes} \dim (A) \geq  h^\Ome_\mu. \EE 
  
  When the  sequence $\EEE$ is composed of wavelet leaders, or of $p$-leaders, the corresponding function spaces are no more Besov spaces, but  alternative families of function spaces, the { \em Oscillation Spaces}, see \cite{JafJMP,BB2}. 

   
 \[ \Hmin = \sup\{ A: (\ref{defab}) \; \mbox{ holds} \}  . \]


The following upper bounds for dimensions are classical for measures, see  \cite{BMP}, and are stated in  the general setting  of dyadic functions in \cite{JARVW}.

\BP \label{sharpuppbou} Let  $\EEE$ be a dyadic function, and let  
\[ J^\Ome_H = \{ x \in \Ome: \; h_\EEE (x) \geq H  \},   \hspace{1mm}  G^\Ome_H =  \{  x \in \Ome : \; h_\EEE (x) \geq H  \}, \] 
  \[  F^\Ome_H = \{ x  \in \Ome : \; \tilde{h}_\EEE (x)   \leq H)\},  \hspace{1mm} K^\Ome_H = \{  x \in \Ome: \; \tilde{h}_\EEE (x)   \leq H)\} .  \]
\begin{itemize}
 \item If $\EEE \in \bsp (\Ome )  $ with $p >0$,  then $\dim (G^\Ome_H) \leq d-sp+Hp$.
\item If $\EEE \in  \btsp (\Ome )  $ with $p >0$,  then $\dim_p (F^\Ome_H) \leq d-sp+Hp$.
\item If $\EEE \in  \bsp (\Ome )  $ with $p <0$,  then $\dim (K^\Ome_H) \leq d-sp+Hp$.
\item If $\EEE \in  \btsp (\Ome )  $ with $p <0$,  then $\dim_p (J^\Ome_H) \leq d-sp+Hp$. 
\end{itemize} 
\EP



 \subsection{Function space  interpretation: Varying regularity  }
 
 Recall that the global scaling function has a function space interpretation in terms of Besov spaces which contain the dyadic function $\EEE$. Similarly, the local scaling function  can be given two functional interpretations; one is local, and  in terms of { \em germ spaces} at  a point, and the second is global, and is in terms of { \em function spaces with varying smoothness}.  We now recall these notions, starting with germ spaces in a general, abstract setting. 
 
 \BD \label{defgerm}  Let $E$ be a Banach space (or a quasi-Banach space)  of distributions  satisfying ${ \mathcal D } \hookrightarrow  E \hookrightarrow { \mathcal D ' }$.  Let $x\in \RR^d$; a  distribution $f$ belongs to $E$ { \bf locally at $x$}  if there exists $\varphi \in D$ such that  $\varphi (x) =1$ in a neighbourhood of $x$ and $f \varphi  \in E$. We also say that $f$ belongs to the  germ space  of $E$ at $x$,  denoted by $E_{x}$. 
 \ED
 
 Let us  draw the relationship between the  local scaling function and  germ spaces: If the $(e_\la)$ are the wavelet coefficients of  a function $f$, then 
 \[ \forall p >0, \hspace{1cm} \eta_f (x, p) = \sup \left\{ s: f \in B^{s/p, \infty}_{p, x} \right\}.   \] 
Note that, in the wavelet case, these  local Besov regularity indices have been investigated by H. Triebel, see Theorem 4  of   \cite{Trie1} where their wavelet characterization is derived, (the reader should be careful that what is referred to as ``pointwise regularity'' in the terminology introduced by H. Triebel  is called here ``local regularity'').

The uniform exponent can also be reformulated  in terms of of H\"older spaces: 
 \[Ê{ \mathcal H}_f (x)  =  \sup \left\{ s: f \in C^{s}_{ x} \right\}.   \] 
 In that case, the function $Ê{ \mathcal H}_f (x)$ is called the{ \em  Local H\"older exponent } of $f$. Its properties have been investigated by S. Seuret and his collaborators, see e.g. \cite{JLVS}.  
 
 Note that the definition of germ spaces can be adapted to the dyadic functions setting.  
 \BD  Let $E$ be a Banach space (or a quasi-Banach space)  defined on  dyadic functions over $\Ome$ ;  a  dyadic function $(e_\la ) $ belongs to $E_x$ if there exists a neighbourhood $\ome$ of $x$ such that the dyadic function $(e_\la ) $  restricted to $\ome$ belongs to $E$. 
 \ED
 
 It is clear that this definition, when restricted to the case of Besov spaces and wavelet coefficients coincides with Definition \ref{defgerm}. 
 \\
 
 We now turn to function spaces with varying smoothness. Such spaces were initially introduced by  Unterberger and  and Bokobza in \cite{UB1,UB2}, followed by many authors (see \cite{Schn} for an extensive review on the subject). A general way to introduce such spaces is to remark that the classical Sobolev spaces $H^{s,p} (\RR^d) $ can be defined by the condition 
 \[ \parallel T(f) \parallel_p  < \infty,  \] 
 where $T$ is the pseudo-differential operator defined by 
 \[   (T f)(x) = \frac{1}{(2\pi)^d}\int_{\RR^d} e^{ix \xi}(1+ | \xi |^2 )^{s/2} \hat{f} (\xi) d \xi  . \]
 This definition  leads to operators with constant order $s$ because the symbol $(1+ | \xi |^2 )^{s/2} $ is independent of $x$. However, one can define more general spaces, with possibly varying order if replacing $(1+ | \xi |^2 )^{s/2}$ by a symbol $\sigma (x, \xi)$. In particular the symbols $(1+ | \xi |^2 )^{a(x)/2}$ will lead to Sobolev   spaces  of varying order $H^{a, p} $ where we can expect that,  if $a$ is  a smooth enough function (say continuous), then the local order of smoothness at $x$  will be $a(x)$. This particular case, and its extensions in the Besov setting, has been studied by H.G. Leopold, followed by J. Schneider,  Besov, H. Triebel, A. Almeida, P. H\"ast\"o,  J. Vyb\'{\i}ral, and several other authors, who gave alternative characterizations of these space in terms of finite differences or Littlewood-Paley decomposition. They also  studied their mutual embeddings (and also in the case where  both the order of smoothness and  the order of integrability $p$ vary)  and their interpolation properties, see  \cite{Schn2,Schn3} fand  references therein, and also \cite{Schn} for an historical account. The reader can also consult \cite{AlHa, Vy} for recent extensions in particular when both the order of smoothness  and the order of integrability $p$ vary. We follow here the  presentation of J. Schneider, since this author obtained Littlewood-Paley characterizations, which are clearly equivalent to the wavelet characterization that we now give. For the sake of simplicity, we assume form now on that the distributions considered are defined on $\RR^d$ and that the wavelet basis used belongs to the Schwartz class (the usual adaptations are standard in the case of functions on a domain, or for wavelets with limited regularity). We additionally assume that the function $a$ is uniformly continuous and satisfies 
 \BE \label{boun}   \exists c, C >0, \; \forall x\in \RR^d, \hspace{1cm}  c \leq a(x) \leq C . \EE
 Then the Besov space $B^{a, q}_p$ (for $p, q \in (0,\infty ]$) can be  characterized by the following wavelet  condition, which is independent of the wavelet basis used.
 
 \BP \label{defbesvar} Let $a$ be a uniformly continuous function satisfying (\ref{boun}), and let $p, q \in (0,\infty ]$.  The  Besov space of varying order  $B^{a, q}_p$ is characterized by the following condition: 
 
 Let $c_\la$ denote the wavelet coefficients of a distribution $f$, and 
 let 
 \[ a_j = \left( 2^{-dj}\sum_{\la \in \La_j} (c_\la 2^{a(\la)  j})^p\right)^{1/p}, \]
 where $a(\la)$ denotes the average of the function $a$ on the cube $\la$; then $f \in B^{a, q}_p$ if 
 $  (a_j ) \in l^q$.  
 \EP

Note that when $p=q=2$ one recovers the Sobolev space $H^{a, 2} $  defined above, and when $a$ is  a constant equal to $s$, then one recovers the standard Besov space $B^{s, q}_p$. Furthermore, the embeddings 
 \[   B^{a, 1}_p \hookrightarrow H^{a, p}  \hookrightarrow B^{a, \infty}_p \]
 yield easy to handle  ``almost characterizations'' of Sobolev spaces of varying order.

The following result, which follows directly from the definition of the local scaling function (Definition \ref{defscalfuncloc}) and the characterization supplied by Proposition \ref{defbesvar}, gives the interpretation of the local scaling function in terms Besov spaces of varying order.

\BP Let $f$ be a distribution defined on $\RR^d$. Then, for $p >0$,  the local wavelet scaling function of $f$ can be recovered by
\[  \forall p >0, \hspace{1cm} \eta_f (p,x) =  p \cdot \sup \{ a: f \in B^{a, \infty}_p \} . \]
\EP

\medskip\medskip
  
{ \bf Acknowledgement: }

This work is supported by the ANR grants AMATIS and MUTATIS and by the LaBeX B\'ezout.  The last two authors thank M. Lapidus for his kind invitation to the fractal geometry session in the AMS conference, March 2012.


\end{document}